\documentclass[preprint,3p,12pt]{elsarticle}

\usepackage{stmaryrd}
\usepackage{bbding}
\usepackage{dcolumn}
\usepackage{graphicx}
\usepackage{amsmath}
\usepackage{amsfonts}
\usepackage{amssymb}
\usepackage{psfrag}
\usepackage{wrapfig}
\usepackage{subfigure}
\usepackage{makeidx}
\usepackage{bm}
\usepackage{epsf}
\usepackage{epsfig}
\usepackage{setspace}
\usepackage{booktabs}
\usepackage{float}
\usepackage{epsfig}
\usepackage{amssymb,amsmath,mathrsfs}
\usepackage{amsopn,wasysym}
\usepackage{subfigure}
\usepackage{color}

\begin{document}

\title{A Third-order Compact Gas-kinetic Scheme on Unstructured Meshes for Compressible Navier-Stokes Solutions}

\author[HKUST1]{Liang Pan}
\ead{panliangjlu@sina.com}
\author[HKUST1,HKUST2]{Kun Xu\corref{cor}}
\ead{makxu@ust.hk} \cortext[cor]{Corresponding author}

\address[HKUST1]{Department of Mathematics, Hong Kong University of Science and Technology, Clear Water Bay, Kowloon, Hong Kong}
\address[HKUST2]{Department of Mechanical and Aerospace Engineering, Hong Kong University of Science and Technology, Clear Water Bay, Kowloon, Hong Kong}

\begin{abstract}
In this paper, for the first time a compact third-order gas-kinetic
scheme is proposed on unstructured meshes for the compressible
viscous flow computations. The possibility to design such a
third-order compact scheme is due to the high-order gas evolution
model, where a time-dependent gas distribution function at a cell
interface not only provides the fluxes across a cell interface, but
also the time evolution of the flow variables at the cell interface
as well. As a result, both cell averaged and cell interface flow
variables can be used for the initial data reconstruction at the
beginning of next time step. A weighted least-square reconstruction
has been used for the construction of a third-order initial
condition. Therefore, a compact third-order gas-kinetic scheme with
the involvement of neighboring cells only can be developed on
unstructured meshes. In comparison with other conventional
high-order schemes, the current method avoids the use of Gaussian
points for the flux integration along a cell interface and the
multi-stage Runge-Kutta time stepping technique. The third-order
compact scheme is numerically stable under CFL condition above
$0.5$. Due to the multidimensional gas-kinetic formulation and the
coupling of inviscid and viscous terms, even with unstructured
meshes the boundary layer solution and the vortex structure can be
accurately captured in the current scheme. At the same time, the
compact scheme can capture strong shocks as well.

\end{abstract}
\begin{keyword}
high-order scheme, gas-kinetic scheme, compact reconstruction,
unstructured mesh, weighted least-square reconstruction.
\end{keyword}
\maketitle

\section{Introduction}

In computational fluid dynamics, the second-order methods are
generally robust and reliable, and they are routinely
employed in the practical calculations. For the same computational
cost, higher-order methods can provide more accurate solutions, but
they are less robust and more complicated. In recent decades, there has been
a continuous interesting and effort on the development of higher-order
schemes. For engineering applications, the construction of higher-order
numerical schemes on unstructured meshes becomes extremely demanding.
Since a gigantic amount of publications have been devoted to the introduction and survey of higher-orders schemes,
the current paper will mainly concentrate on the construction of the
third-order compact gas-kinetic scheme on unstructured meshes.

The gas-kinetic scheme (GKS) has been developed systematically for
the compressible flow computations \cite{GKS-Xu1,
GKS-Xu2,GKS-Kumar,GKS-Jiang}. An evolution process from  kinetic to
hydrodynamic scales has been constructed for the flux evaluation.
The kinetic effect through particle free transport contributes to
the capturing of the shock wave, and the hydrodynamic effect plays a
dominant role for the resolved viscous and heat conducting
solutions. In other words, the highly non-equilibrium of the gas
distribution function in the discontinuous region provides a
physically consistent mechanism for the construction of a numerical
shock structure. In this sense, the GKS is close to the methodology
of artificial dissipation approach, but with different dissipative
mechanism. In smooth flow region, the hydrodynamic scale physics
corresponding to the multi-dimensional central difference
discretization captures the accurate viscous solutions. Due to the
coupling of inviscid and viscous terms in the kinetic formulation,
theoretically there is no difficulty for GKS to capture NS solutions
in any structure or unstructured mesh. With the discretization of
particle velocity space, a unified gas-kinetic scheme (UGKS) has
been developed for the flow study in entire flow regime from
rarefied to continuum ones \cite{UGKS-Xu,UGKS-Luc,UGKS-Guo}.

Recently, with the incorporation of high-order initial data
reconstruction, a higher-order gas-kinetic schemes has been proposed
in \cite{GKS-high1, GKS-high2,GKS-high3}. The flux evaluation in the
scheme is based on the time evolution of flow variables from an
initial piece-wise discontinuous polynomials (parabola) around a
cell interface, where higher-order spatial and temporal derivatives
of a gas distribution function are coupled nonlinearly. The whole
curves of discontinuous flow distributions around a cell interface
interact through particle transport and collision in the
determination of the flux function. Besides the evaluation of the
time-dependent flux function across a cell interface, the
higher-order gas evolution model also provides an accurate
time-dependent solution of flow variables at a cell interface as
well. Thus, it is feasible to develop a compact scheme with the
consideration of time evolution of both cell averaged  and cell
interface flow variables. A compact third-order gas-kinetic scheme
is proposed for the compressible Euler and Navier-Stokes equations
on structure meshes with WENO-type reconstruction \cite{GKS-high4}.
However, this reconstruction technique is difficult to be used on
unstructured meshes. Therefore, on the unstructured meshes, a
weighted least-square reconstruction will be used  in this paper. To
the third-order accuracy, a quadratic distribution for the flow
variables inside each cell needs to be determined. Based on the cell
averaged and cell interface values of neighboring cells only, an
over-determined linear system is formed. With the least-square
solution for the system, the whole  flow distribution can be fully
determined. The shock detector can be also used as well to switch
between higher-order (3rd) and lower order (2nd) reconstructions in
different regions. In comparison with traditional schemes, the
Gaussian points for the flux evaluation along the cell interface and
the multi-stage Runge-Kutta technique are avoided in the current
compact method. At the same time, the current third-order compact
scheme is stable under the CFL condition  $\mbox{CFL} \simeq 0.5$.

This paper is organized as follows. In Section 2, the finite volume
scheme on the unstructured mesh and third-order GKS are
introduced. In section 3, the compact reconstruction on the triangular mesh
is presented, and the techniques can be applied to rectangular mesh as well.
Section 4 includes numerical
examples to validate the current algorithm. The last section is the
conclusion.

\section{Finite volume gas-kinetic scheme}
\subsection{Finite volume scheme}
The two-dimensional gas-kinetic BGK equation can be written as \cite{BGK-1},
\begin{equation}\label{bgk}
f_t+\textbf{u}\cdot\nabla f=\frac{g-f}{\tau},
\end{equation}
where $f$ is the gas distribution function, $g$ is the corresponding
equilibrium state, and $\tau$ is the collision time. The collision
term satisfies the compatibility condition
\begin{equation}\label{compatibility}
\int \frac{g-f}{\tau}\varphi d\Xi=0,
\end{equation}
where $\varphi=(1,u,v,\displaystyle \frac{1}{2}(u^2+v^2+\xi^2))$,
$d\Xi=dudvd\xi^1...d\xi^{K}$, $K$ is the number of internal freedom,
i.e.  $K=(4-2\gamma)/(\gamma-1)$ for two-dimensional flows, and
$\gamma$ is the specific heat ratio.

Based on the Chapman-Enskog expansion of the BGK model, the Euler
and Navier-Stokes, Burnett, and Super-Burnett equations can be
derived \cite{BGK-3, GKS-Xu1}. In the smooth region, the gas
distribution function can be expanded as
\begin{align*}
f=g-\tau D_{\textbf{u}}g+\tau D_{\textbf{u}}(\tau
D_{\textbf{u}})g-\tau D_{\textbf{u}}[\tau D_{\textbf{u}}(\tau
D_{\textbf{u}})g]+...,
\end{align*}
where $D_{\textbf{u}}={\partial}/{\partial
t}+\textbf{u}\cdot \nabla$. By truncating different orders of
$\tau$, the corresponding macroscopic equations can be derived. For
the Euler equations, the zeroth order truncation is taken, i.e.
$f=g$. For the Navier-Stokes equations, the first order truncation
is
\begin{align}\label{ns}
f=g-\tau (ug_x+vg_y+g_t).
\end{align}
Based on the higher order truncations, the Burnett and super-Burnett
eqautions can be obtained.

In the computation, the computational volumes are simply triangles.
For a control volume $\Omega_i$, its boundary is given by three
line segments
\begin{equation*}
\partial\Omega_i=\bigcup_m\Gamma_{im}.
\end{equation*}
Thus, taking moments of the kinetic equation Eq.\eqref{bgk} and
integrating with respect to time and space, the finite volume scheme
can be expressed as
\begin{align}\label{finite}
W_{i}^{n+1}=W_{i}^{n}-\frac{1}{|\Omega_i|}\int_{t^n}^{t^{n+1}}\sum_mF_{im}(t)dt,
\end{align}
where $W=(\rho,\rho U,\rho V,\rho E)$ are the conservative
variables, $F_{im}(t)=(F_{\rho},F_{\rho u},F_{\rho v},F_{E})$ are
the fluxes across the cell interface $\Gamma_{im}$ in the global
coordinate, which is defined as
\begin{align}\label{flux1}
F_{im}(t)=\int_{\Gamma_{im}}(\int\varphi f(x,y,t,u,v,\xi) \textbf{u}
\cdot \textbf{n}du dvd\xi)ds.
\end{align}
where $\textbf{n}=(\cos\theta,\sin\theta)$ is the outer normal
direction of the cell interface $\Gamma_{im}$, and the tangential
direction is denoted as $\textbf{t}=(-\sin\theta,\cos\theta)$.
Eq.\eqref{finite} is valid in any scale if the interface flux is
properly defined, which is beyond the validity of the Navier-Stokes
equations.

According to the coordinate transformation, the local coordinate for
the cell interface $\Gamma_{im}$ is expressed as
$(\widetilde{x},\widetilde{y})=(0, \widetilde{y})$, where
$\widetilde{y}\in[-d, d]$ and $d=|\Gamma_{im}|/2$, and the
velocities in the local coordinate are given by
\begin{align}\label{uu1}
\begin{cases}
\widetilde{u}=u\cos\theta+v\sin\theta,\\
\widetilde{v}=-u\sin\theta+v\cos\theta .
\end{cases}
\end{align}
For the gas distribution function in the local coordinate,
$\widetilde{f}(\widetilde{x},\widetilde{y},t,\widetilde{u},\widetilde{v},\xi)=f(x,y,t,u,v,\xi)$
and $dudv=d\widetilde{u}d\widetilde{v}$, then the line integral for
the gas distribution function over the cell interface $\Gamma_{im}$
can be transformed as
\begin{align}\label{flux3} \int_{\Gamma_{im}}\int\varphi
f(x,y,t,u,v,\xi) \textbf{u} \cdot \textbf{n}du dvd\xi
ds=\int_{-d}^{d}\int\varphi\widetilde{f}(0,\widetilde{y},t,\widetilde{u},\widetilde{v},\xi)\widetilde{u}d\widetilde{u}d\widetilde{v}d\xi
d\widetilde{y}.
\end{align}
Thus, in the computation, the numerical fluxes in the local
coordinate
$\widetilde{F}(t)=(F_{\widetilde{\rho}},F_{\widetilde{m}},F_{\widetilde{n}},F_{\widetilde{E}})$
are obtained first by taking moments of the gas
distribution function in the local coordinate
\begin{align}\label{flux2}
\widetilde{F}(t)=\int_{-d}^{d}\int\widetilde{u}\widetilde{\varphi}
\widetilde{f}(0,\widetilde{y},t,\widetilde{u},\widetilde{v},\xi)d\widetilde{u}d\widetilde{v}d\xi
d\widetilde{y},
\end{align}
where
$\widetilde{\varphi}=(1,\widetilde{u},\widetilde{v},\displaystyle\frac{1}{2}(\widetilde{u}^2+\widetilde{v}^2+\xi^2))$.
According to Eq.\eqref{uu1} and Eq.\eqref{flux3}, the fluxes in the
global coordinate can be expressed as a combination of the fluxes in
the local coordinate
\begin{align}\label{flux-1}
\begin{cases}
F_{\rho}=F_{\widetilde{\rho}},\\
F_{m}=F_{\widetilde{m}}\cos\theta-F_{\widetilde{n}}\sin\theta,\\
F_{n}=F_{\widetilde{m}}\sin\theta+F_{\widetilde{n}}\cos\theta,\\
F_{E}=F_{\widetilde{E}}.
\end{cases}
\end{align}
With the above numerical fluxes at the cell interface, the flow
variables inside each control volume can be updated according to
Eq.\eqref{finite}.

\subsection{Gas-kinetic flux solver}
In this section, the numerical flux will be presented in the local
coordinate. For simplicity, all notations
with tilde will be omitted here after.

In order to simulate the NS solutions, we need to model the interface flux function.
For the distribution function at a cell interface, the integral solution of BGK
equation Eq.\eqref{bgk} at the cell interface in the local
coordinate can be written as
\begin{equation}\label{integral1}
f(0,y,t,u,v,\xi)=\frac{1}{\tau}\int_0^t g(x',y',t',u,v,\xi)e^{-(t-t')/\tau}dt'\\
+e^{-t/\tau}f_0(-ut,y-vt,u,v,\xi),
\end{equation}
where $x=0$ is the location of the cell interface, $x=x'+u(t-t')$
and $y=y'+v(t-t')$ are the trajectory of particles, $f_0$ is the
initial gas distribution function, and $g$ is the corresponding
equilibrium state. The target equations to be solved depend on the modeling of the initial condition $f_0$ term.

To construct a multidimensional third-order gas-kinetic solver, the
following notations are introduced firstly
\begin{align*}
a_1=&(\partial g/\partial x)/g, a_2=(\partial g/\partial y)/g,
A=(\partial g/\partial t)/g, B=(\partial A /\partial t),\\
d_{11}&=(\partial a_1/\partial x), d_{12}=(\partial a_1/\partial
y)=(\partial a_2/\partial x), d_{22}=(\partial a_2/\partial y),
\\
&b_{1}=(\partial a_1/\partial t)=(\partial A/\partial x),
b_{2}=(\partial a_2/\partial t)=(\partial A/\partial y),
\end{align*}
where $g$ is an equilibrium state. The dependence of these
coefficients on particle velocity can be expanded as the following
form \cite{GKS-Xu2}
\begin{align*}
a_1=a_{11}+a_{12}u+&a_{13}v+a_{14}\displaystyle
\frac{1}{2}(u^2+v^2+\xi^2),\\
&...\\
B=B_{1}+B_{2}u+&B_{3}v+B_{4}\displaystyle
\frac{1}{2}(u^2+v^2+\xi^2).
\end{align*}

For the kinetic part of the integral solution Eq.\eqref{integral1},
the gas distribution function can be constructed as
\begin{equation}\label{f0}
f_0=f_0^l(x,y,u,v)H(x)+f_0^r(x,y,u,v)(1-H(x)),
\end{equation}
where $H(x)$ is the Heaviside function,  $f_0^l$ and $f_0^r$ are the
initial gas distribution functions on both sides of a cell
interface, which have one to one correspondence with the initially
reconstructed polynomials of macroscopic flow variables on both
sides of the cell interface. To construct a third-order scheme, the
Taylor expansion for the gas distribution function in space and time
at $(x,y)=(0,0)$ is expressed as
\begin{align*}
f_0^k(x,y)=f_G^k(0,0)&+\frac{\partial f_G^k}{\partial
x}x+\frac{\partial f_G^k}{\partial y}y+\frac{1}{2}\frac{\partial^2
f_G^k}{\partial x^2}x^2+\frac{\partial^2 f_G^k}{\partial x\partial
y}xy+\frac{1}{2}\frac{\partial^2 f_G^k}{\partial y^2}y^2,\nonumber
\end{align*}
where $k=l,r$. For the Euler equations, $f_{G}^k=g_k$ and the
kinetic part of Eq.\eqref{integral1} can be obtained. For the
Navier-Stokes equations, according to Eq.\eqref{ns} and the
notations introduced above, the distribution function is
\begin{align*}
f_{G}^k=g_k-\tau(a_{1k}u+a_{2k}v+A_k)g_k,
\end{align*}
where $g_l,g_r$ are the equilibrium states corresponding to the
macroscopic variables $W_l, W_r$ given by the reconstruction
procedure at both sides of cell interface. Thus, the corresponding
kinetic part of Eq.\eqref{integral1} can be written as
\begin{align}
&e^{-t/\tau}f_0^k(-ut,y-vt,u,v)\nonumber\\
=&C_7g_k[1-\tau(a_{1k}u+a_{2k}v+A_k)]\nonumber\\
+&C_8g_k[a_{1k}u-\tau((a_{1k}^2+d_{11k})u^2+(a_{1k}a_{2k}+d_{12k})uv+(A_ka_{1k}+b_{1k})u)]\nonumber\\
+&C_8g_k[a_{2k}v-\tau((a_{1k}a_{2k}+d_{12k})uv+(a_{2k}^2+d_{22k})v^2+(A_ka_{2k}+b_{2k})v)]\nonumber\\
+&C_7g_k[a_{2k}-\tau((a_{1k}a_{2k}+d_{12k})u+(a_{2k}^2+d_{22k})v+(A_ka_{2k}+b_{2k}))]y\nonumber\\
+&\frac{1}{2}C_7g_k[(a_{1k}^2+d_{11k})(-ut)^2+2(a_{1k}a_{2k}+d_{12k})(-ut)(y-vt)+(a_{2k}^2+d_{22k})(y-vt)^2],\label{dis2}
\end{align}
where $g_{k}$ are the equilibrium states at both sides of the cell
interface, and the coefficients $a_{1k},...,A_k$ are defined
according to the expansion of $g_{k}$.

After determining the kinetic part $f_0$, the equilibrium state $g$
in the integral solution Eq.\eqref{integral1} can be constructed
 as follows
\begin{align}\label{equli}
g=g_0+\frac{\partial g_0}{\partial x}x+&\frac{\partial g_0}{\partial
y}y+\frac{\partial g_0}{\partial t}t+\frac{1}{2}\frac{\partial^2
g_0}{\partial x^2}x^2+\frac{\partial^2 g_0}{\partial x\partial
y}xy+\frac{1}{2}\frac{\partial^2 g_0}{\partial
y^2}y^2\nonumber\\
&+\frac{1}{2}\frac{\partial^2 g_0}{\partial t^2}t^2+\frac{\partial^2
g_0}{\partial x\partial t}xt+\frac{\partial^2 g_0}{\partial
y\partial t}yt,
\end{align}
where $g_{0}$ is the equilibrium state located at interface, which
can be determined through the compatibility condition
Eq.\eqref{compatibility}
\begin{align}\label{compatibility2}
\int\psi g_{0}d\Xi=W_0=\int_{u>0}\psi g_{l}d\Xi+\int_{u<0}\psi
g_{r}d\Xi.
\end{align}
Based on Taylor expansion for the equilibrium state
Eq.\eqref{equli}, the hydrodynamic part in Eq.\eqref{integral1} can
be written as
\begin{align}\label{dis1}
\frac{1}{\tau}\int_0^t
g&(x',y',t',u,v)e^{-(t-t')/\tau}dt'\nonumber\\
=&C_1g_0+C_2g_0\overline{a}_1u+C_2g_0\overline{a}_2v+C_1g_0\overline{a}_2y+C_3g_0\overline{A}\nonumber\\
+&\frac{1}{2}C_4g_0(\overline{a}_1^2+\overline{d}_{11})u^2+C_6g_0(\overline{A}\overline{a}_1+\overline{b}_{1})u+\frac{1}{2}C_5g_0(\overline{A}^2+\overline{B})\nonumber\\
+&\frac{1}{2}C_1g_0(\overline{a}_2^2+\overline{d}_{22})y^2+C_2g_0(\overline{a}_2^2+\overline{d}_{22})vy+\frac{1}{2}C_4g_0(\overline{a}_2^2+\overline{d}_{22})v^2\nonumber \\
+&C_2g_0(\overline{a}_1\overline{a}_2+\overline{d}_{12})uy+C_4g_0(\overline{a}_1\overline{a}_2+\overline{d}_{12})uv\nonumber\\
+&C_3g_0(\overline{A}\overline{a}_2+\overline{b}_{2})y+C_6g_0(\overline{A}\overline{a}_2+\overline{b}_{2})v,
\end{align}
where the coefficients
$\overline{a}_1,\overline{a}_2,...,\overline{A},\overline{B}$ are
defined from the expansion of the equilibrium state $g_0$. The
coefficients $C_i, i=1,...,8$ in Eq.\eqref{dis1} and Eq.\eqref{dis2}
are given by
\begin{align*}
C_1=1-&e^{-t/\tau}, C_2=(t+\tau)e^{-t/\tau}-\tau, C_3=t-\tau+\tau e^{-t/\tau},C_4=-(t^2+2t\tau)e^{-t/\tau},\\
&C_5=t^2-2t\tau,C_6=-t\tau(1+e^{-t/\tau}),C_7=e^{-t/\tau},C_8=-te^{-t/\tau}.
\end{align*}
Substituting Eq.\eqref{dis1} and Eq.\eqref{dis2} into the integral
solution Eq.\eqref{integral1}, the gas distribution function at the
cell interface can be obtained. The superscripts or subscripts of
the coefficients $a_1, a_2,...,A, B$ in Eq.\eqref{dis2} and
Eq.\eqref{dis1} are omitted for simplicity and they are determined
by the spatial derivatives of macroscopic flow variables and the
compatibility condition \cite{GKS-high2} as follows
\begin{align}\label{var-fun}
\begin{cases}
\displaystyle\langle a_1\rangle =\frac{\partial W}{\partial x},
\langle a_2\rangle =\frac{\partial W}{\partial y},  \langle
A+a_1u+a_2v \rangle=0,\\ \displaystyle\langle a_1
^2+d_{11}\rangle=\frac{\partial^2 W}{\partial x^2}, \langle a_2
^2+d_{22}\rangle=\frac{\partial^2 W}{\partial y^2}, \langle
a_1a_2+d_{12}\rangle=\frac{\partial^2
W}{\partial x\partial y},\\
\displaystyle\langle(a_1 ^2+d_{11})u+(a_1a_2+d_{12})v+(Aa_1+b_1)\rangle=0,\\
\displaystyle\langle(a_1a_2+d_{12})u+(a_2 ^2+d_{22})v+(Aa_2+b_2)\rangle=0,\\
\displaystyle\langle(Aa_1+b_1)u+(Aa_2+b_2)v+(A^2+B)\rangle=0,
\end{cases}
\end{align}
where $<...>$ are the moments of gas distribution function, and
defined by
\begin{align*}
<...>=\int g(...)\psi d\Xi.
\end{align*}
In the following section, with the reconstruction procedure, the
conservative value $W_{l}, W_{r}$ and $W_{0}$ at the center of cell
interface corresponding to the equilibrium $g_l, g_r, g_0$ and the
derivatives in Eq.\eqref{var-fun} will be presented.

\section{Compact reconstruction}
\begin{figure}[!h]
\centering
\includegraphics[width=0.4\textwidth]{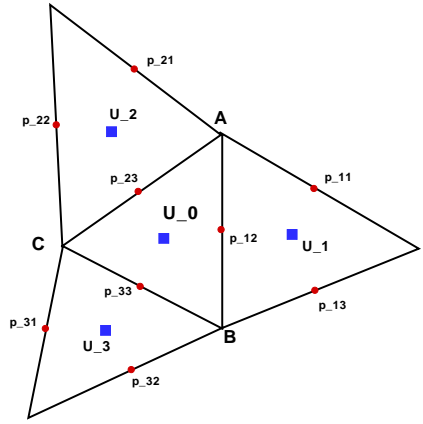}
\caption{\label{compact-s} The stencil of a compact reconstruction
for triangle $\Omega_0=\vartriangle_{ABC}$. The blue squares are the
cell averaged values and the red circles are point values at the
center of cell interface.}
\end{figure}

This paper focuses on the high-order compact finite volume scheme.
In the finite volume type schemes, to achieve higher-order accuracy,
a reconstruction for the flow variables with high-order polynomials
inside each cell is needed as the initial condition at the beginning
of each time step. For the higher-order reconstruction, a large
number of stencils is usually needed to determine all degrees of
freedom through the WENO or least square techniques
\cite{un-ENO,un-WENO2,un-WENO3,k-exact-1,k-exact-3}. In this
section, the reconstruction will be done  for the unstructured mesh
with a compact stencil, which is shown in Fig.{\ref{compact-s}}. For
simplicity, the whole reconstruction procedure is performed in a
local coordinate $(x,y)$ relative to a cell interface, such as  AB
in Fig.{\ref{compact-s}}, which is consistent with the evaluation of
a time-dependent gas distribution function at the cell interface.

In the gas-kinetic scheme, besides the numerical
fluxes, the macroscopic pointwise values at a cell interface in the
local coordinate can be obtained by taking moments of the gas
distribution function,
\begin{align}\label{point}
W(t,y)&=\int \varphi f(x_{i+1/2},y,t,u,v,\xi)dudvd\xi.
\end{align}
As shown in the last section, the whole curve of the polynomial of
the macroscopic variables will participate the flow evolution, and
the spatial and temporal derivatives of the gas distribution
function are coupled nonlinearly. This point-wise value at the cell
interface Eq.\eqref{point} is a solution of the evolution model,
which can be used in the reconstruction stage at the beginning of
next time step. Thus, in the following subsections, a third-order
compact reconstruction will be presented for the unstructured mesh,
in which the pointwise values at the cell interface and the cell
averaged values shown in Fig.{\ref{compact-s}} are used in the
reconstruction.

The macroscopic variables for reconstruction is denoted by $U$. For
the smooth flow, the conservative variables $W$ will be directly
used for reconstruction, i.e. $U=W$. For the flow with
discontinuity, in order to eliminate the spurious oscillation and improve the
stability of the scheme, the compact reconstruction is based on the
characteristic variables. Denote $F(W)=(\rho U, \rho U^2+p, \rho
UV,U(\rho E+p))$ in the local coordinate. The Jacobian matrix
$\partial F/\partial W$ can be diagnoalized by the right eigenmatrix
$R$, and the characteristic variables is defined as $U=R^{-1}W$. For
a cell interface, $R$ is the right eigenmatrix for $\partial
F/\partial W^*$ and $W^*$ is the averaged conservative value from
both side of cell interface. To the third order accuracy, the
expansion of the macroscopic variable $U$ inside the cell $\Omega_0$ can
be expressed as
\begin{align}\label{expansion}
U(x,y)=U_0&+U_{x}((x-x_0)-\widehat{x}_0)+U_{y}((y-y_0)-\widehat{y}_0)+\frac{1}{2}U_{xx}((x-x_0)^2-\widehat{x_0^2})\nonumber\\
+&U_{xy}((x-x_0)(y-y_0)-\widehat{x_0y_0})+\frac{1}{2}U_{yy}((y-y_0)^2-\widehat{y_0^2}),
\end{align}
where $(x_0,y_0)$ is the barycenter of $\Omega_0$, $U_0$ is the cell
averaged value for $U(x,y)$, and
\begin{align*}
\widehat{x^my^n}=\frac{1}{|\Omega_0|}\int_{\Omega_0}(x-x_0)^n(y-y_0)^mdV.
\end{align*}
The cell averaged value for the base function over the triangle
$\Omega_i$ is denoted as
\begin{align}\label{cell-va}
\widehat{x^my^n}_i=\frac{1}{|\Omega_i|}\int_{\Omega_i}((x-x_0)^n(y-y_0)^m-\widehat{x^my^n})dV.
\end{align}
and the point-wise value for the base function at the point
$p_{ij}=(x_{ij},y_{ij})$ is denoted as
\begin{align}\label{point-va}
x^my^n_{ij}=((x_{ij}-x_0)^n(y_{ij}-y_0)^m-\widehat{x^my^n}).
\end{align}

\subsection{Initial data reconstruction}
In this subsection, the weighted least-square reconstruction will be
presented for the initial data reconstruction. As shown in
Fig.\ref{compact-s}, three cell averaged values $U_i, i=1,2,3$ (blue
square) form the neighboring cells and nine point-wise values
$U_{ij}, i,j=1,2,3$  (red circle) from the cell interface will be
used in the weighted least square reconstruction.

For the third order expansion, with the definition of the cell
averaged and point-wise values for the base function
Eq.\eqref{cell-va} and Eq.\eqref{point-va}, we have
\begin{align}\label{third-re1}
U_x\widehat{x}_i+U_y\widehat{y}_i+\frac{1}{2}U_{xx}\widehat{x^2}_i+U_{xy}\widehat{xy}_i+\frac{1}{2}U_{yy}\widehat{y^2}_i=U_i-U_0,
\end{align}
where $U_i$ is the cell averaged value for the neighboring triangle
$\Omega_i$, $i=1,2,3$. For the nine cell interface points $p_{ij}$,
$i,j=1,2,3$, we have
\begin{align}\label{third-re2}
U_xx_{ij}+U_yy_{ij}+\frac{1}{2}U_{xx}x^2_{ij}+U_{xy}xy_{ij}+\frac{1}{2}U_{yy}y^2_{ij}=U_{ij}-U_0,
\end{align}
where $U_{ij}$ is the point-wise value of $U(x,y)$ at the point
$p_{ij}$.

To solve the corresponding derivatives for $U(x,y)$,
Eq.\eqref{third-re1} and Eq.\eqref{third-re2} can be written into
an over-determined linear system
\begin{align}\label{over-determined}
\displaystyle\left(\begin{array}{ccccc}
\widehat{x}_1&\widehat{y}_1&\frac{1}{2}\widehat{x^2}_1&\widehat{xy}_1&\frac{1}{2}\widehat{y^2}_1\\
~&~&...&~&~\\
\widehat{x}_3&\widehat{y}_3&\frac{1}{2}\widehat{x^2}_3&\widehat{xy}_3&\frac{1}{2}\widehat{y^2}_3\\
x_{11}&y_{11}&\frac{1}{2}x^2_{11}&xy_{11}&\frac{1}{2}y^2_{11}\\
~&~&...&~&~\\
x_{33}&y_{33}&\frac{1}{2}x^2_{33}&xy_{33}&\frac{1}{2}y^2_{33}
\end{array}\right)\cdot \left(\begin{array}{c}
U_x\\
U_y\\
U_{xx}\\
U_{xy}\\
U_{yy}
\end{array}
\right)= \left(\begin{array}{c}
U_1-U_0 \\
...\\
U_3-U_0 \\
U_{{11}}-U_0 \\
...\\
U_{{33}}-U_0
\end{array}\right).
\end{align}
Denote $dU=(U_x, U_y, U_{xx},U_{xy}, U_{yy})^T$, $\Delta
U=(U_1-U_0,...,U_3-U_0,U_{{11}}-U_0,...,U_{{33}}-U_0)^T$, the above
linear system is expressed as the matrix form
\begin{align*}
DdU=\Delta U.
\end{align*}
where $D$ is the coefficient matrix corresponding to
Eq.\eqref{over-determined}.

To deal with the discontinuity, a diagonal matrix $W$ is introduced
as the simple weight functions
\begin{align*}
w_{i}=\frac{1}{(s_i^2+\epsilon)},~~~
w_{ij}=\frac{1}{(s_{ij}^2+\epsilon)}
\end{align*}
where $\displaystyle
s_i=\frac{U_i-U_0}{|\mathbf{x}_i-\mathbf{x}_0|},
s_{ij}=\frac{U_{ij}-U_0}{|\mathbf{x}_i-\mathbf{x}_0|}$,
$i,j=1,...,3$, and $\epsilon=10^{-6}$. The derivatives $dU$ can be
obtained by solving the linear system
\begin{eqnarray*}
D^TWDdU=D^TW\Delta U.
\end{eqnarray*}

Generally, for most cases with Mach number $Ma<2$, the weight
function is enough to deal with the discontinuity. However, for strong discontinuity, the shock detection \cite{Shock-detection}
technique is used in the current scheme. Analogous to the analysis
of KXRCF detector \cite{Shock-detection}, for the third-order
scheme, it is easy to distinguish the smooth region from the region
near discontinuities as follows
\begin{align*}
U_{i}(x_{i})-U_{j\rightarrow i}(x_i)=\begin{cases}
O(h^{3})~~\text{in smooth region,}\\
O(h)~~~\text{near discontinuity,}
\end{cases}
\end{align*}
where the index $i$ refers $\Delta_{ABC}$ and the index $j$ refers
$\Delta_{ABC'}$, $U_{i}(x_{i})$ is the interpolated value at the
center of $\Delta_{ABC}$ and $U_{j\rightarrow i}(x_{i})$ is the
value at the center of $\Delta_{ABC}$ extrapolated from
$\Delta_{ABC'}$. In the computation, the "trouble cell" is detected
according the following criterion
\begin{align*}
\max(\|U_{i}(x_{i})-U_{j\rightarrow
i}(x_i)\|,\|U_{j}(x_{j})-U_{i\rightarrow j}(x_j)\|)\geq
C\sqrt{S_{\Delta_{ABC}}+S_{\Delta_{ABC'}}}\sim O(h),
\end{align*}
where $S$ is the area of the triangle, $C$ is a problem dependent
coefficient, and $C=5$ is used in the computation. In those detected
"trouble cell", the second order scheme with limiters are used.
The above choice of weight functions may not be optimal and further study is needed.

\begin{figure}[!h]
\centering
\includegraphics[width=0.5\textwidth]{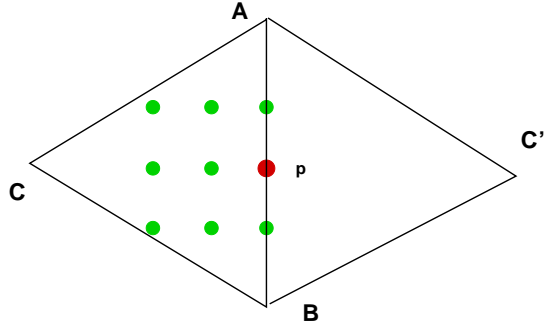}
\caption{\label{compact-s1} The stencil of the compact
reconstruction for triangle $\Omega_0=\vartriangle_{ABC}$ for the
characteristic variables.  The coordinate of these points $p_{ij},
i,j=1,2,3$ is $((i-3)d,(j-2)d)$, where $d=d_{AB}/4$.}
\end{figure}

With the derivatives $dU=(U_x, U_y, U_{xx},U_{xy}, U_{yy})^T$, the
whole flow distribution in the cell $\Delta_{ABC}$ in Fig.\ref{compact-s1} can be obtained.
For the smooth flow, no
special treatment is needed. With $W=U$, the interpolated value
$W_l$ and the derivatives $dW_l$ can be fully obtained in the cell
$\Delta_{ABC}$ . Similarly, the interpolated value $W_r$ and the
derivatives $dW_r$ in the cell $\Delta_{ABC'}$ can be obtained as well.

For the flow with discontinuity, the characteristic variables are
reconstructed in the cell $\Delta_{ABC}$.  With the derivatives
$dU=(U_x, U_y, U_{xx},U_{xy}, U_{yy})^T$, the interpolated value $U$
at the points in Fig.\ref{compact-s1} can be obtained. By the
inverse projection, the conservative variables $W=RU$, where $R$ is
the right eigenmatrix. Based on these point-wise values and their
central difference, $W_l$ and $dW_l$ can be obtained. Similarly, the
interpolated value $W_r$ and the derivatives $dW_r$ in the cell
$\Delta_{ABC'}$ can be also obtained.

\subsection{Reconstruction for equilibrium part}
In this subsection, the reconstruction for the equilibrium part will
be presented. This reconstruction will be based on the conservative
variables $W$. To the third-order accuracy, the Taylor expansion
corresponding to equilibrium part at the center point of a cell
interface is expressed as
\begin{align}\label{expansion2}
&\overline{W}(x,y)=W_0+\overline{W}_{x}(x-x_p)+\overline{W}_{y}(y-y_p)\nonumber \\
+\frac{1}{2}\overline{W}_{xx}&(x-x_p)^2+\overline{W}_{xy}(x-x_p)(y-y_p)+\frac{1}{2}\overline{W}_{yy}(y-y_p)^2,
\end{align}
where $W_0$ is the conservative variable at the center point of cell
interface $AB$ based on the compatibility condition
Eq.\eqref{compatibility2}, and $\overline{W}_{x}, ...,
\overline{W}_{yy}$ are corresponding derivatives.

\begin{figure}[!h]
\centering
\includegraphics[width=0.35\textwidth]{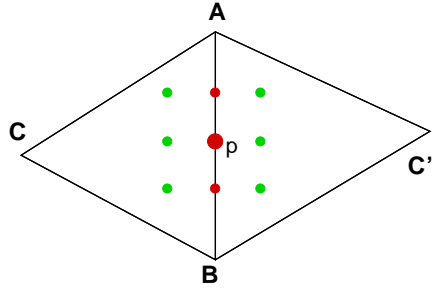}
\caption{\label{schematic-2} The stencil for the equilibrium part in
the local coordinate. The coordinate of these points $p_{ij},
i,j=1,2,3$ is $((i-2)d,(j-2)d)$, where $d=d_{AB}/4$. }
\end{figure}

As shown in Fig.\ref{schematic-2}, with the reconstructed
polynomials in $\vartriangle_{ABC}$ and $\vartriangle_{ABC'}$, the
point values at those points can be determined, which has been
obtained in the last subsection. Especially, we can get the
point values at the interface (red) points at both sides of $AB$. By
the compatibility condition Eq.\eqref{compatibility}, the
reconstructed conservative variables at the cell interface can be
determined. The derivatives $\overline{W}_{x}, ...,
\overline{W}_{yy}$ can be obtained by the central difference of
these point-wise values.

%For each cell interface $\Gamma_{ik}$, the conservative variables
%$W$ are projected into the local characteristic space, which is
%denoted by $U$. With the reconstructed curved $U(x,y)$ and
%$\overline{U}(x,y)$, some point-wise reconstructed conservative
%variables can be obtained by the inverse of characteristic
%projection. The derivatives in Eq.\eqref{var-fun} can be given by
%the numerical derivatives based on these point-wise values, and the
%conservative variables corresponding to the equilibrium $g_l, r_r$
%and $g_0$ in Eq.\eqref{dis2} and Eq.\eqref{dis2} can also obtained
%by the inverse of projection for $U_{pl}$, $U_{pr}$ and $U_{p0}$.
%Thus, the compact reconstruction has been presented.

\begin{figure}[!h]
\centering
\includegraphics[width=0.35\textwidth]{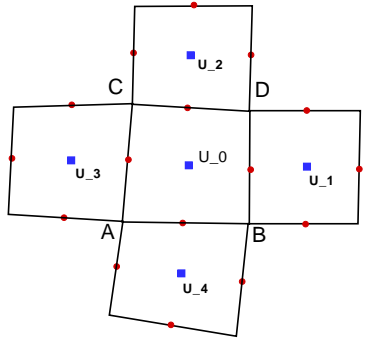}
\caption{\label{schematic-r} The stencil for the rectangular mesh.
The red circles represent the point-wise value and the blue squares
are the cell averaged values.}
\end{figure}

\subsubsection{Extension to rectangular mesh}
For the rectangular mesh, the stencils are given in
Fig.\ref{schematic-r}. To reconstruct the polynomial for the
rectangular $\square_{ABCD}$, the cell averaged values $U_i$,
$i=1,2,3,4$ and point-wise values  $p_{ij}$, $i,j=1,2,3,4$ at the
cell interfaces can be used. Similar to the triangular case,
we have the following matrix form for the over-determined linear
system
\begin{align*}
D_2dU=\Delta U.
\end{align*}
where $dU=(U_x, U_y, U_{xx},U_{xy}, U_{yy})^T$, $\Delta
U=(U_1-U_0,...,U_4-U_0,U_{{11}}-U_0,...,U_{{44}}-U_0)^T$. $D_2$ is
the coefficient matrix and expressed as
\begin{align*}
D_2=\displaystyle\left(\begin{array}{ccccc}
\widehat{x}_1&\widehat{y}_1&\frac{1}{2}\widehat{x^2}_1&\widehat{xy}_1&\frac{1}{2}\widehat{y^2}_1\\
~&~&...&~&~\\
\widehat{x}_4&\widehat{y}_4&\frac{1}{2}\widehat{x^2}_4&\widehat{xy}_4&\frac{1}{2}\widehat{y^2}_4\\
x_{11}&y_{11}&\frac{1}{2}x^2_{11}&xy_{11}&\frac{1}{2}y^2_{11}\\
~&~&...&~&~\\
x_{14}&y_{44}&\frac{1}{2}x^2_{44}&xy_{44}&\frac{1}{2}y^2_{44}
\end{array}\right)
\end{align*}
By introducing the weight diagonal matrix $W$, the derivative $dU$
can be also obtained by solving the following linear system
\begin{eqnarray*}
D_2^TWD_2dU=D_2^TW\Delta U.
\end{eqnarray*}
The limiting process is also used for the flow with large
discontinuity. In some cases of the numerical tests, the solutions from the compact scheme with rectangular mesh will be presented as well.

\section{Numerical tests}
In this section, numerical tests for both inviscid flow and viscous
flow will be presented to validate the compact scheme. For the
inviscid flow, the collision time $\tau$ takes
\begin{align*}
\tau=\epsilon \Delta t+C\displaystyle|\frac{p_l-p_r}{p_l+p_r}|\Delta
t,
\end{align*}
where $\varepsilon=0.05$ and $C=1$. For the viscous flow, we have
\begin{align*}
\tau=\frac{\mu}{p}+C\displaystyle|\frac{p_l-p_r}{p_l+p_r}|\Delta t,
\end{align*}
where $p_l$ and $p_r$ denotes the pressure on the left and right
sides of the cell interface, $\mu$ is the viscous coefficient,  $p$
is the pressure at the cell interface and $C=1$. In the smooth flow
regions, it will reduce to $\tau=\mu/p$. The ratio of specific heats
takes $\gamma=1.4$. $\Delta t$ is the time step which is determined
according to the CFL condition. In the numerical tests, the CFL
number takes a value of $0.35$, even though the scheme works as well
with a large CFL number. The value of $0.35$ is already more than
two times of the time step used for the compact third-order DG
method.

\subsection{Accuracy test}
The numerical order of the compact gas-kinetic scheme is tested in
comparison with the analytical solutions of the Euler equations. The
isotropic vortex propagation problem is presented to validate the
accuracy for the solution of inviscid flow. The computational domain is taken to
be $[0, 1.5]\times[0, 1.5]$. The free upstream is $(\rho, u, v, p) =
(1.21, 0, 0, 1)$, and a small vortex is obtained through a
perturbation on the mean flow with the velocity $(u, v)$,
temperature $T=p/\rho$, and entropy $S=\ln(p/\rho^\gamma)$.
The perturbation is expressed as
\begin{align*}
&(\delta u,\delta v)=\kappa\eta e^{\mu(1-\eta^2)}(\sin\theta,-\cos\theta),\\
&\delta
T=-\frac{(\gamma-1)\kappa^2}{4\mu\gamma}e^{2\mu(1-\eta^2)},\delta
S=0,
\end{align*}
where $\eta=r/r_c$, $r=\sqrt{(x-x_c)^2+(y-y_c)^2}$, $(x_c,
y_c)=(0.75, 0.75)$, $\kappa=0.3$, $\mu=0.204$, and $r_c=0.05$. In the
computation, the unstructured meshes with mesh size $h=1/30, 1/50,
1/100$ and $1/200$ are used, and the $L^\infty$ errors and orders at
$t=1$ are presented in Table.\ref{tab1}, which shows a third-order
accuracy of the current compact scheme.

\begin{table}[!h]
\begin{center}
\def\temptablewidth{0.5\textwidth}
{\rule{\temptablewidth}{0.5pt}}
\begin{tabular*}{\temptablewidth}{@{\extracolsep{\fill}}c|cc}
mesh & $L^\infty$ norm & order ~   \\
\hline
1/30  & 3.2460690E-03 &  ~~       \\
1/50  & 7.3230267E-04 &  2.914901  \\
1/100 & 9.2029572E-05 &  2.992271  \\
1/200 & 1.1801720E-05 &  2.963100
\end{tabular*}
{\rule{\temptablewidth}{0.5pt}}
\end{center}
\vspace{-4mm} \caption{\label{tab1} Accuracy test for the isotropic
vortex problem.}
\end{table}

\begin{figure}[!h]
\centering
\includegraphics[width=0.6\textwidth]{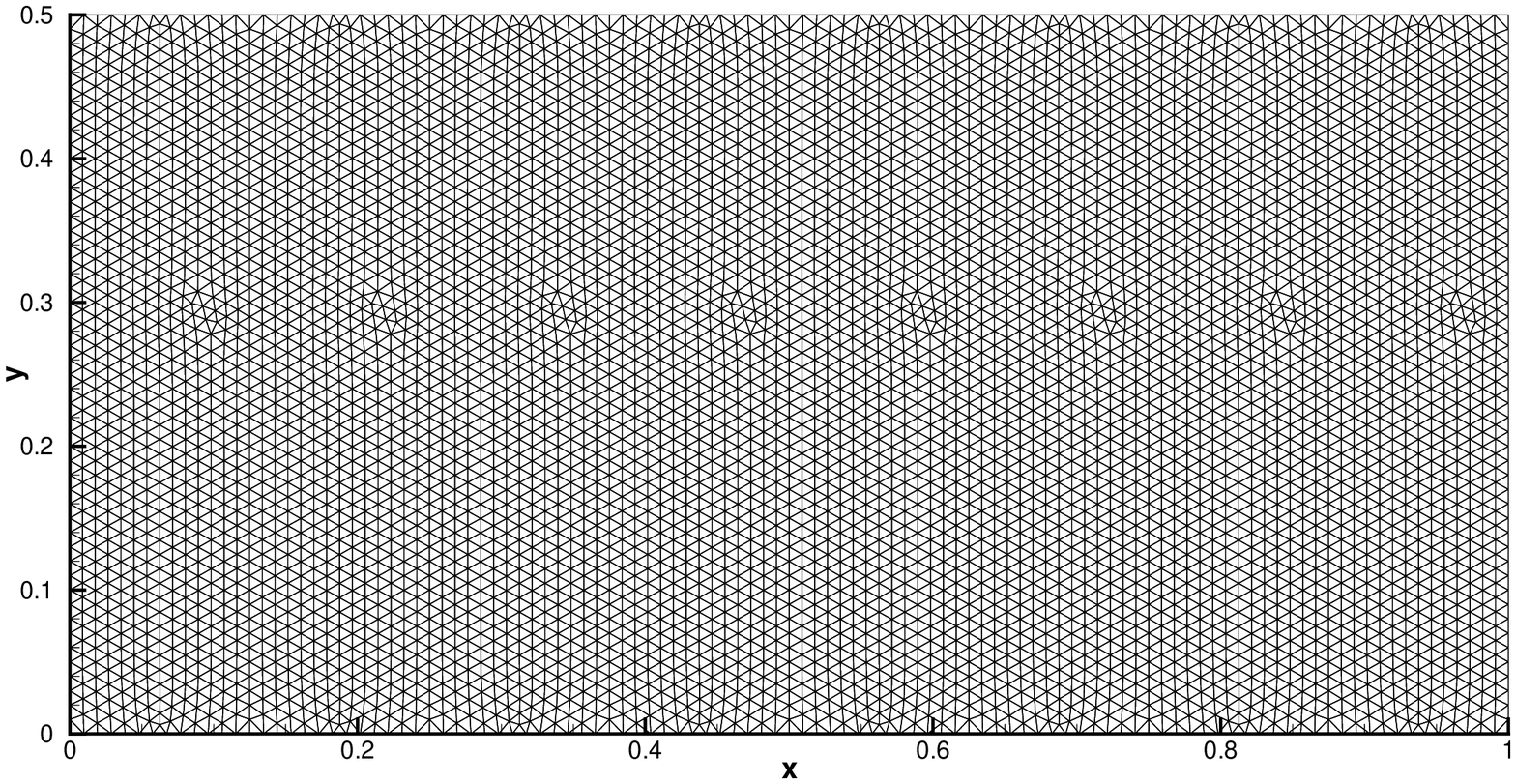}
\caption{\label{riemann-1} 1D Riemann problem: the mesh for the 1D
Riemann problem.}
\includegraphics[width=0.4\textwidth]{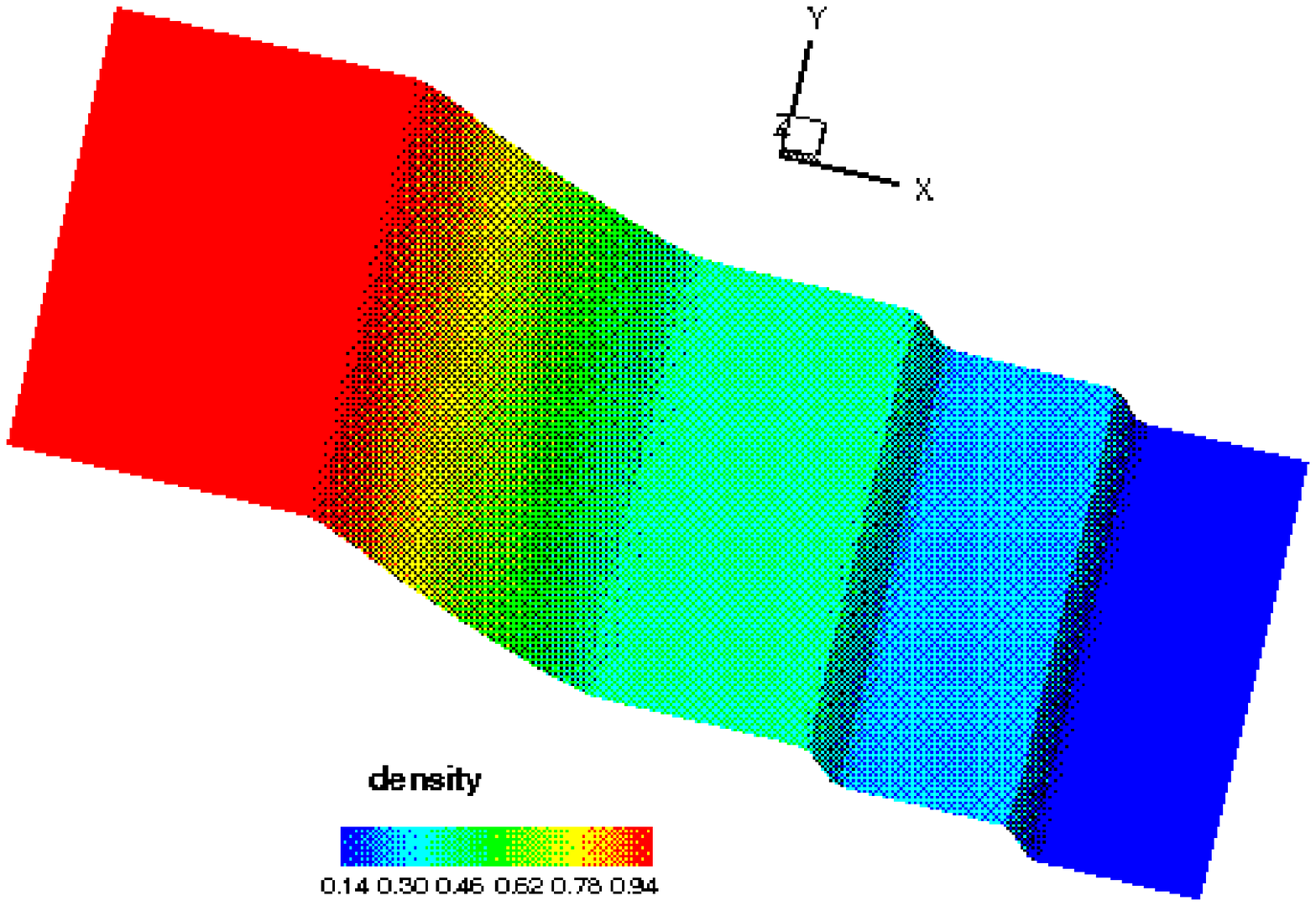}
\includegraphics[width=0.4\textwidth]{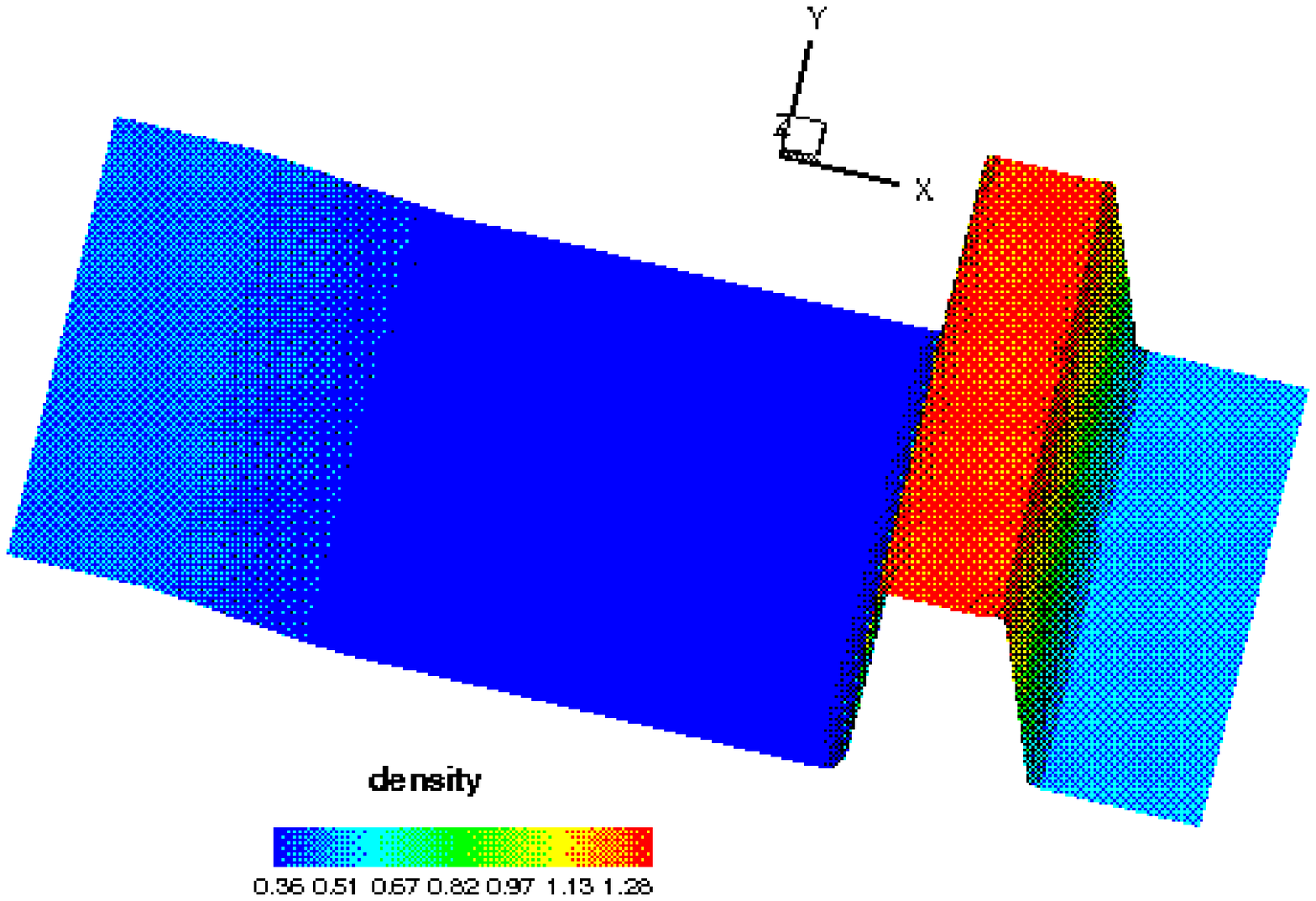}
\caption{\label{riemann-3} 1D Riemann problem: the 3d density
distribution for the Sod problem (left) and Lax problem (right) in
the computational domain.}
\end{figure}
\begin{figure}[!h]
\centering
\includegraphics[width=0.44\textwidth]{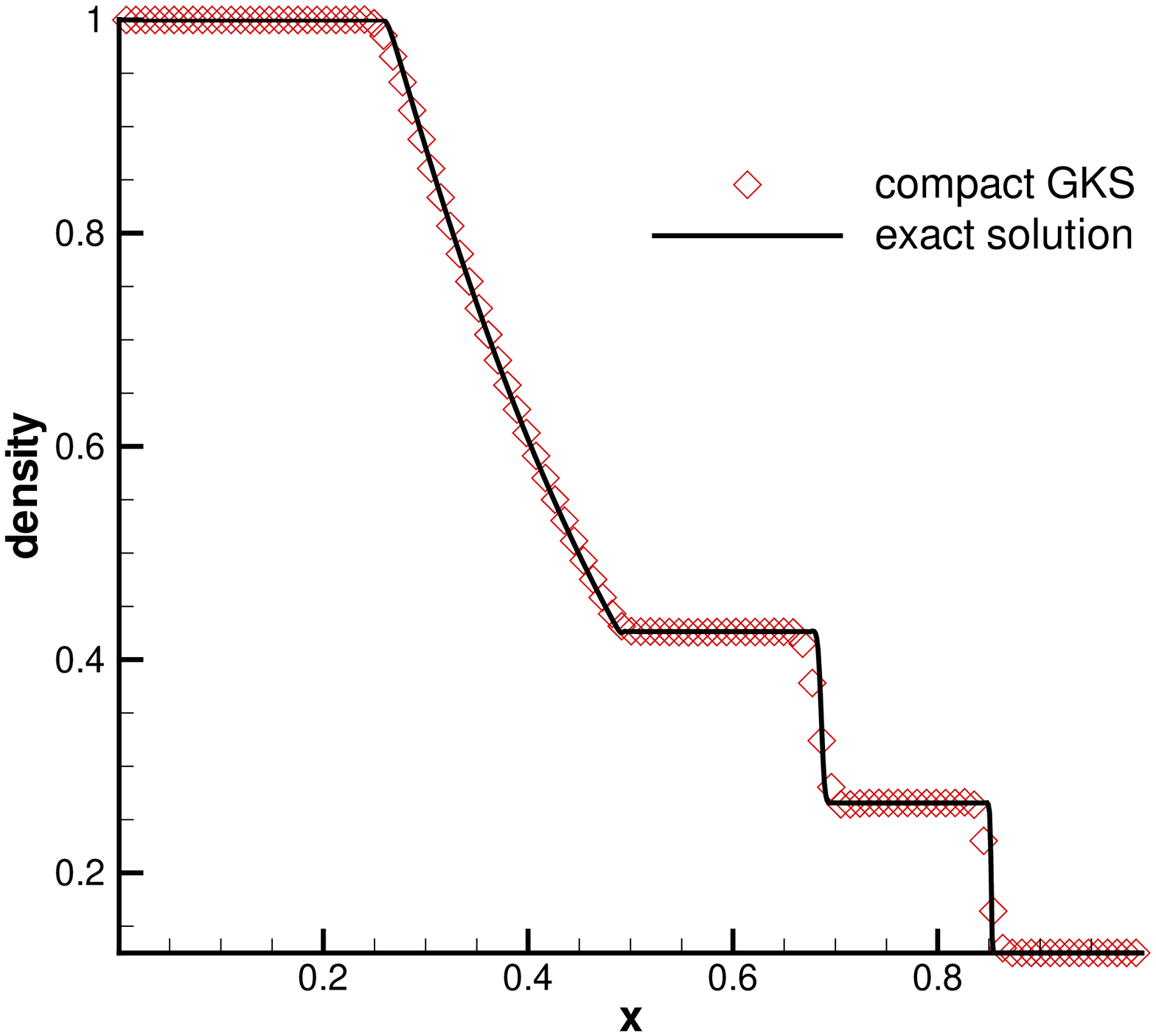}\includegraphics[width=0.44\textwidth]{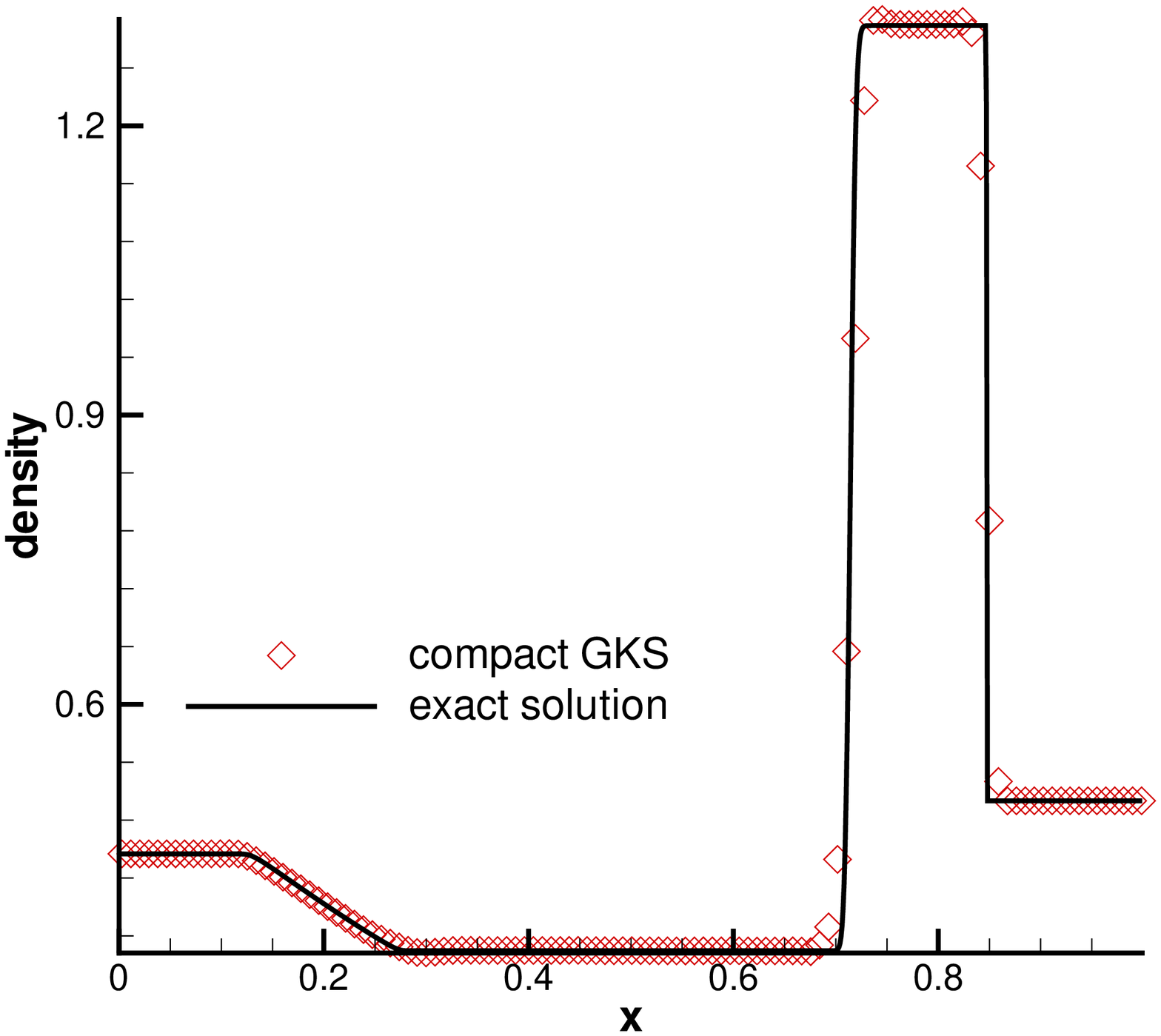}
\includegraphics[width=0.44\textwidth]{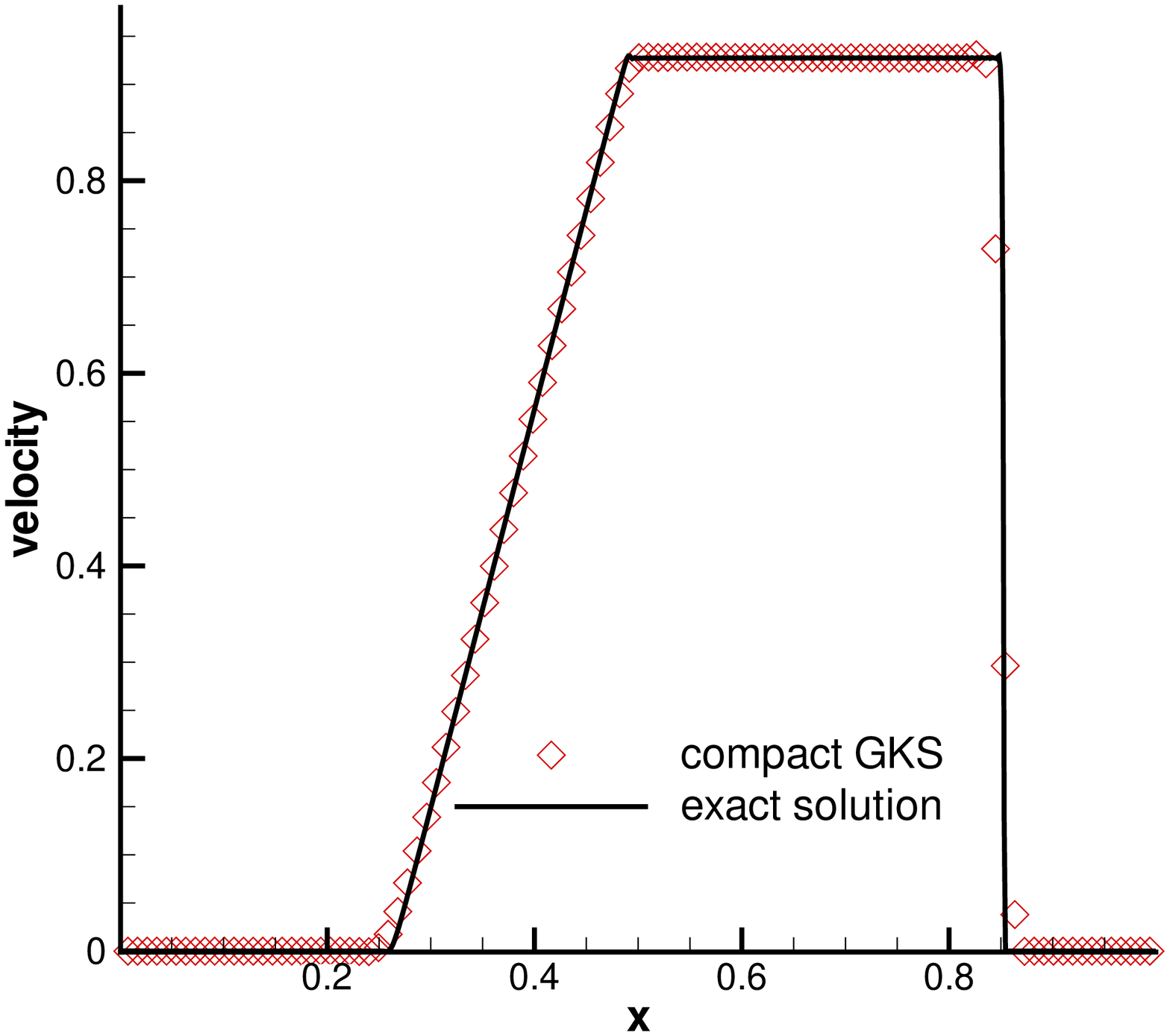}\includegraphics[width=0.44\textwidth]{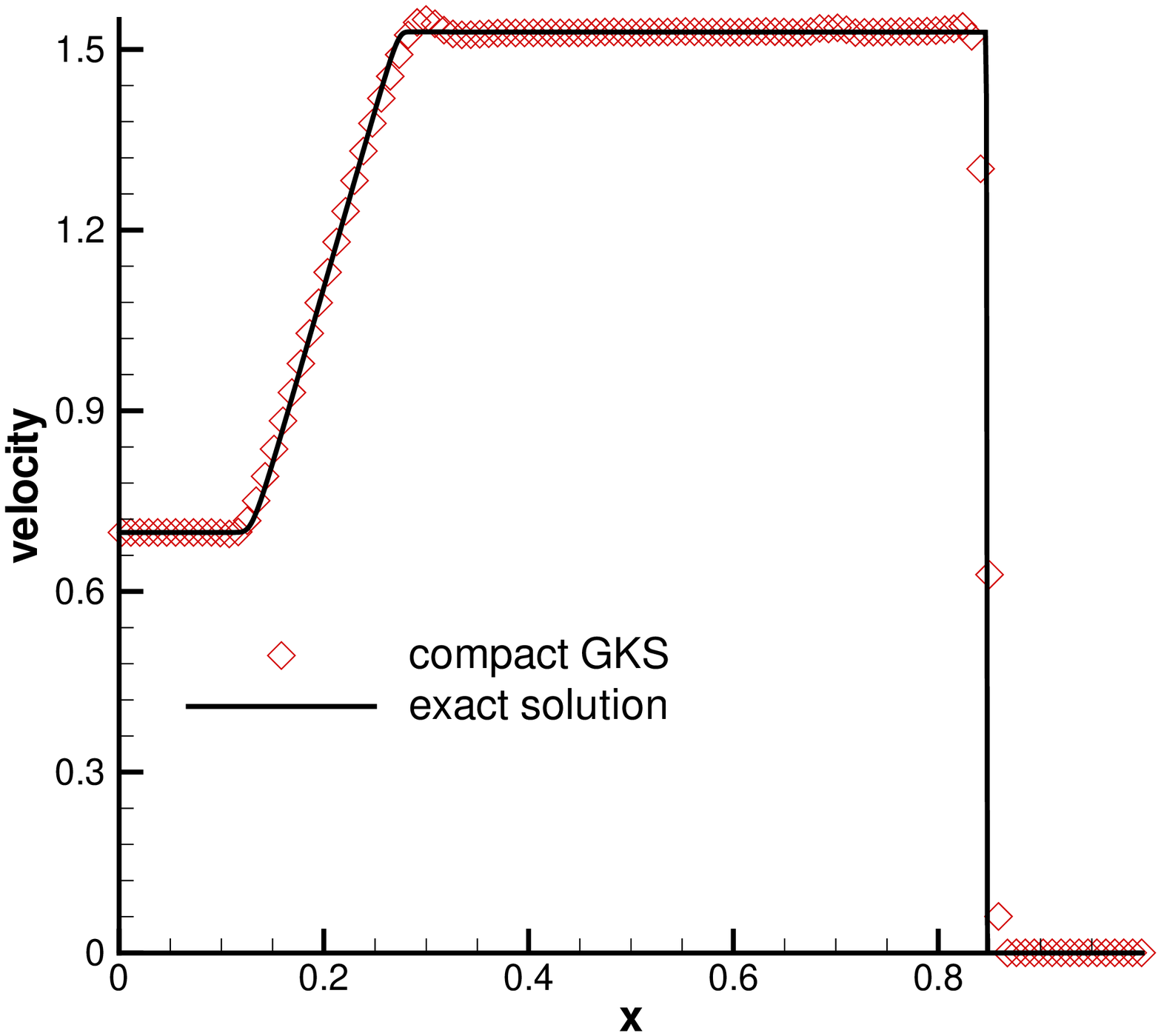}
\includegraphics[width=0.44\textwidth]{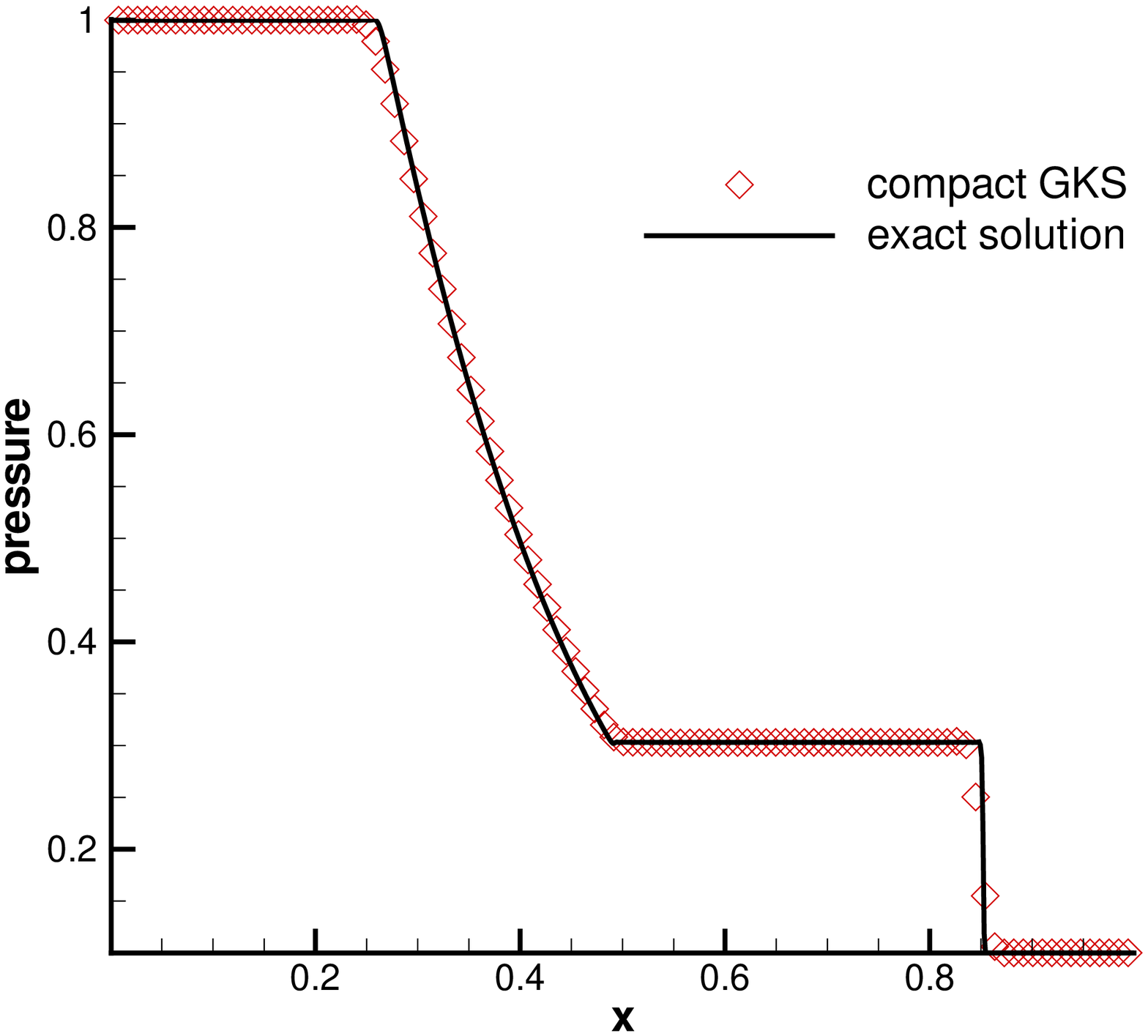}\includegraphics[width=0.44\textwidth]{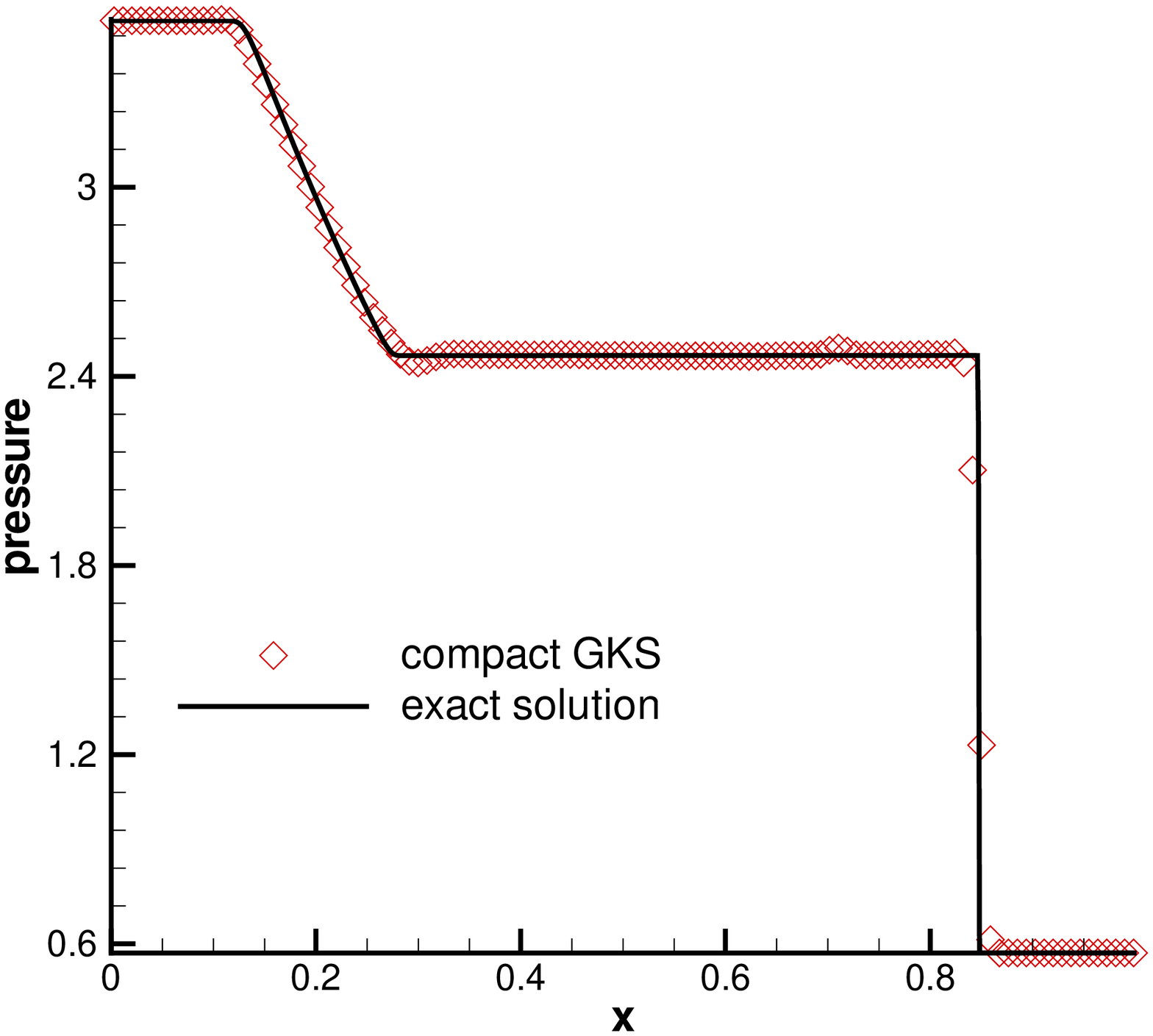}
\caption{\label{riemann-2}1D Riemann problem: Sod problem (left):
the density, velocity, and pressure distributions at t=0.2, and Lax
problem (right): the density, velocity, and pressure distributions at
$t=0.14$, where the mesh size is $h=1/100$.}
\end{figure}

\subsection{One dimensional Riemann problem}
In this case, two one-dimensional Riemann problems are tested to
verify the capability in capturing the wave configurations. The mesh
is presented in Fig.\ref{riemann-1}, where the computational domain is
$[0,1]\times[0,0.5]$, and mesh size is around $h=0.01$. The first one is Sod
problem, and the initial condition is given by
\begin{equation*}
(\rho,u,p)=\left\{\begin{aligned}
&(1, 0, 1), 0<x<0.5,\\
&(0.125,0,0.1),  0.5<x<1.
\end{aligned} \right.
\end{equation*}
The second one is the Lax problem, and the initial condition is
given as follows
\begin{equation*}
(\rho,u,p)=\left\{\begin{aligned}
&(0.445,0.698,3.528), 0\leq x<0.5,\\
&(0.5,0,0.571),  0.5\leq x\leq 1.
\end{aligned} \right.
\end{equation*}
To compare with the exact solution, $100$ points were extracted at $y=0.25$ for the Sod problem at $t=0.2$ and, for the
Lax problem at $t=0.14$. The density, velocity, and pressure
distributions for the exact solutions and numerical results are
presented in Fig.\ref{riemann-2}, where the numerical results agree
well with the exact solutions. The three dimensional density
distributions for the two cases are given in Fig.\ref{riemann-3}. In
this case, the weighted least square reconstruction can deal with the
discontinuity well, and the shock detection technique is not needed.

\subsection{Flow impinging on a blunt body}
\begin{figure}[!h]
\centering
\includegraphics[height=0.166\textwidth]{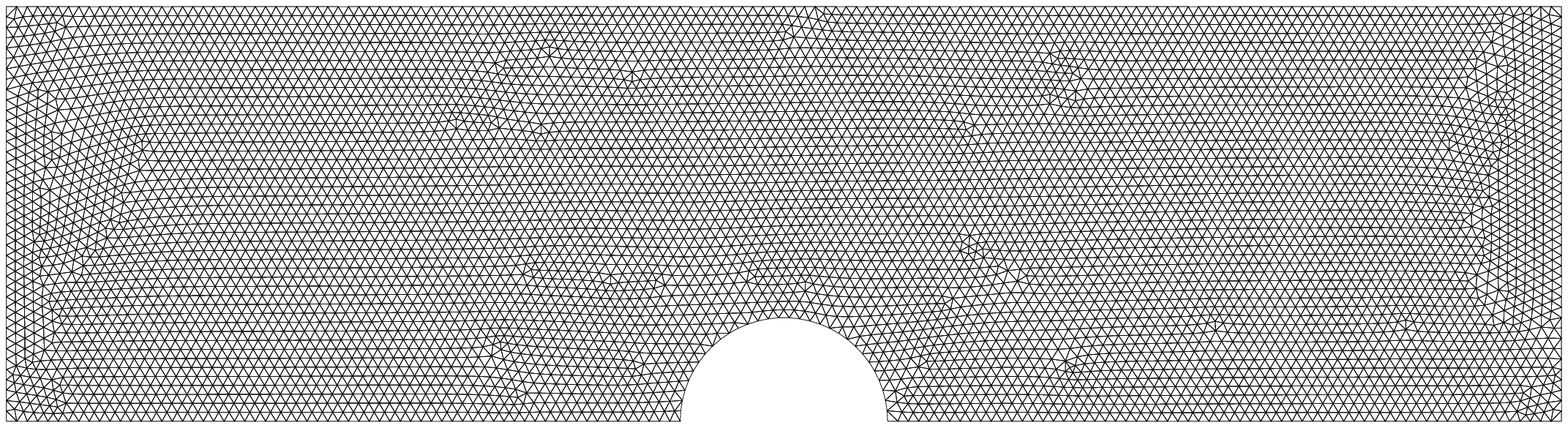}\\
\includegraphics[height=0.166\textwidth]{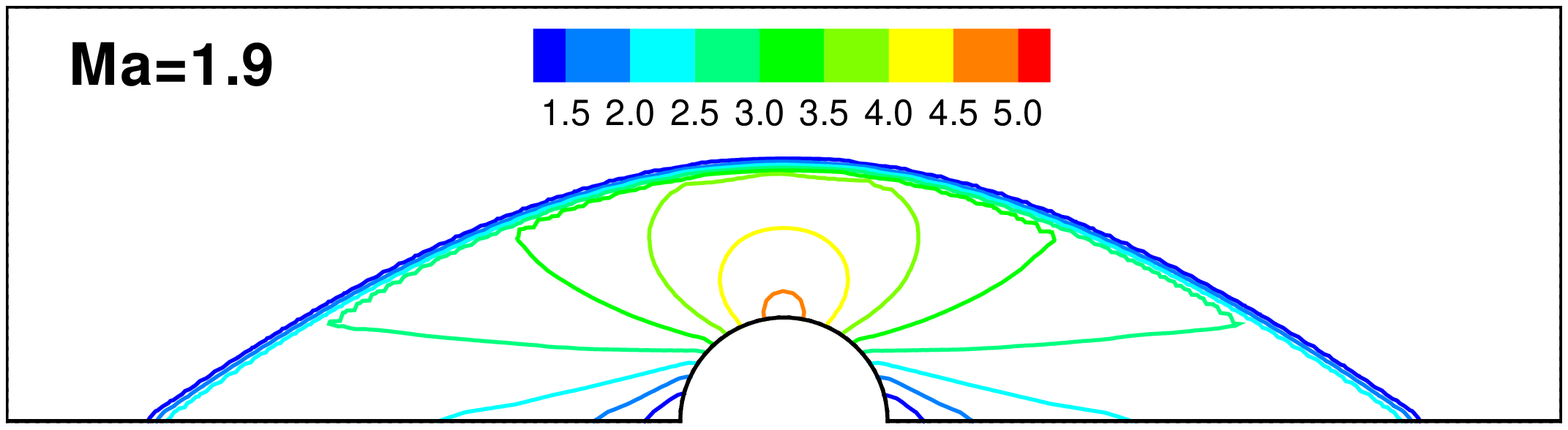}
\caption{\label{cylinder1}Flow impinging on a blunt body: the mesh
and pressure distribution at $Ma=1.9$. }
\includegraphics[height=0.166\textwidth]{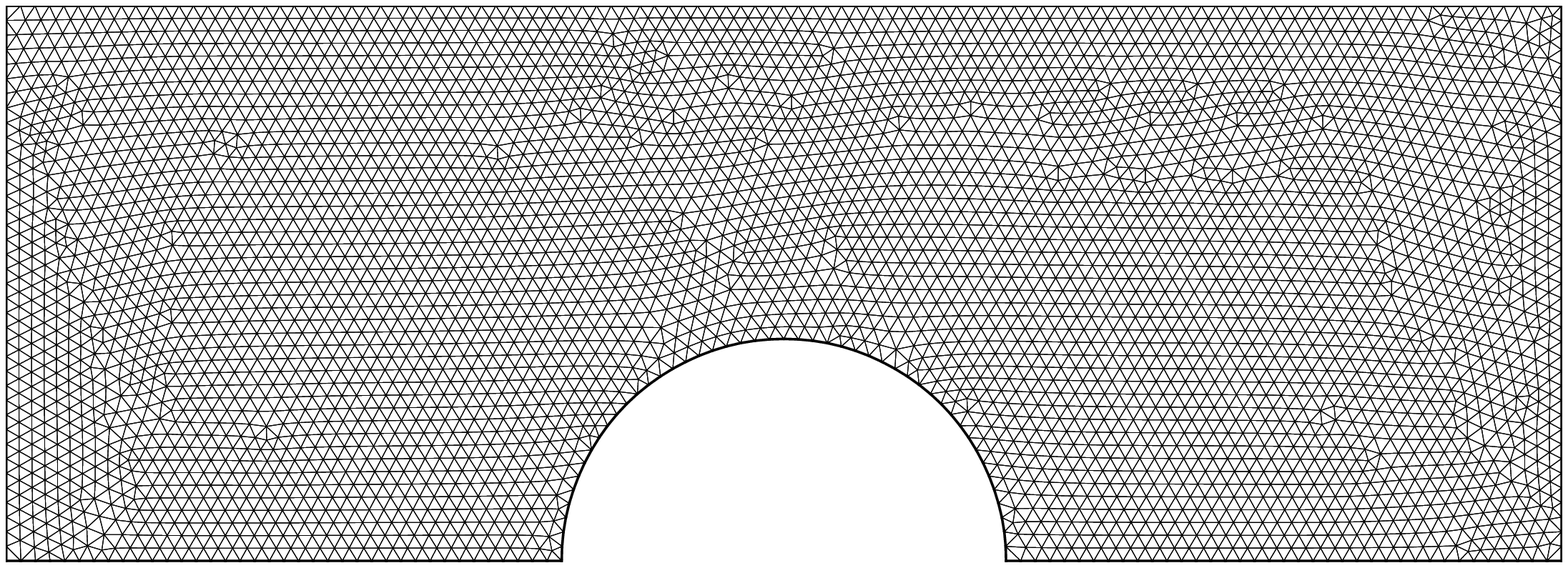}\\
\includegraphics[height=0.166\textwidth]{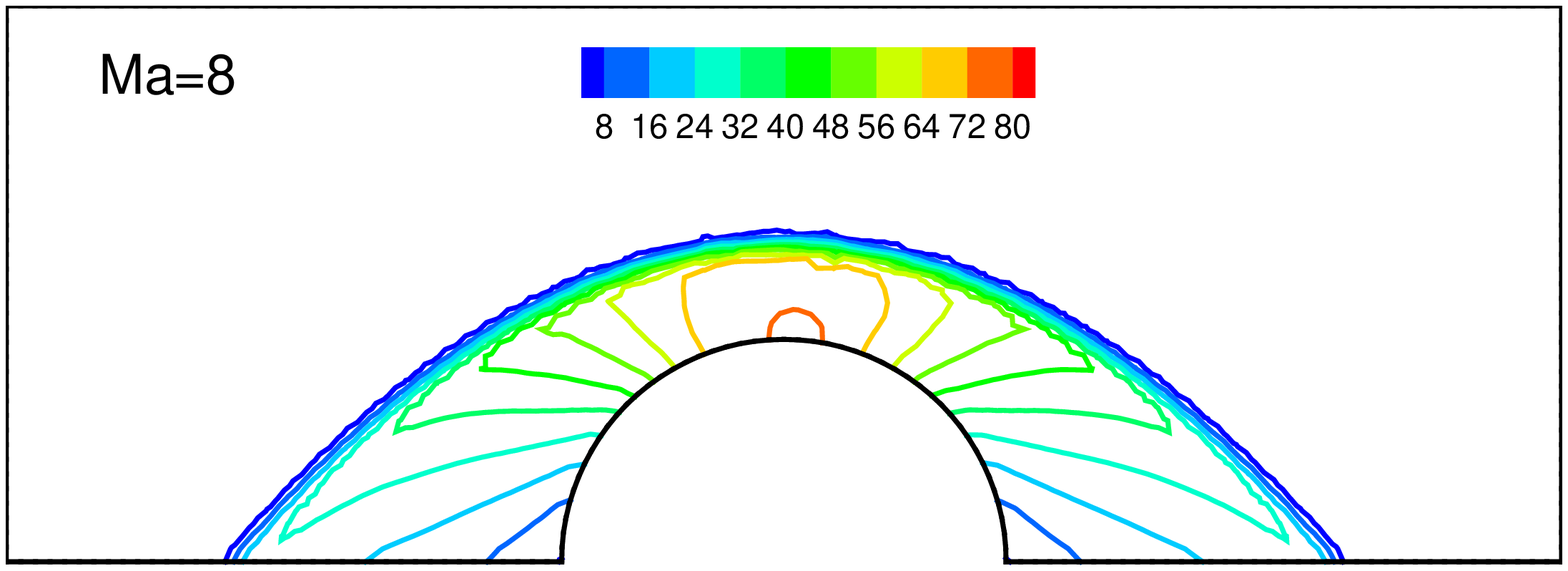}
\caption{\label{cylinder2}Flow impinging on a blunt body: the mesh
and pressure distribution at $Ma=8$. }
\end{figure}

In this case, the inviscid hypersonic flows impinging on a unit
cylinder are tested to validate robustness of the current scheme.
This problem is initialized by the flow moving towards a cylinder
with different Mach numbers. The Euler  boundary condition is
imposed on the surface of cylinder, and outflow boundary condition
on the right boundary.  As mentioned in the reconstruction part,
the weighted least square reconstruction is able to deal with the
discontinuities at a Mach number $Ma<2$. In this case, the flow with
$Ma=1.9$ is tested without the detection of "trouble cell". The mesh
and the pressure distribution for this case are also given in
Fig.\ref{cylinder1}, with mesh size $h=1/15$, where the flow structure
can be captured nicely in front of the cylinder. However, with a
high  Mach number, the weighted least square reconstruction is no
longer able to capture strong discontinuities, and the shock
detection technique is used to identify the trouble cells, where a second-order reconstruction is used in these cells.
For the flow with $Ma=8$, the mesh and
the pressure distribution are shown in
Fig.\ref{cylinder2} with mesh size $h=1/15$. This test shows that the
current scheme can capture the flow structure nicely in front of the
cylinder and the carbuncle phenomenon does not appear
\cite{Case-Pandolfi}.

\subsection{Shock vortex interaction}
\begin{figure}[!h]
\centering
\includegraphics[width=0.4\textwidth]{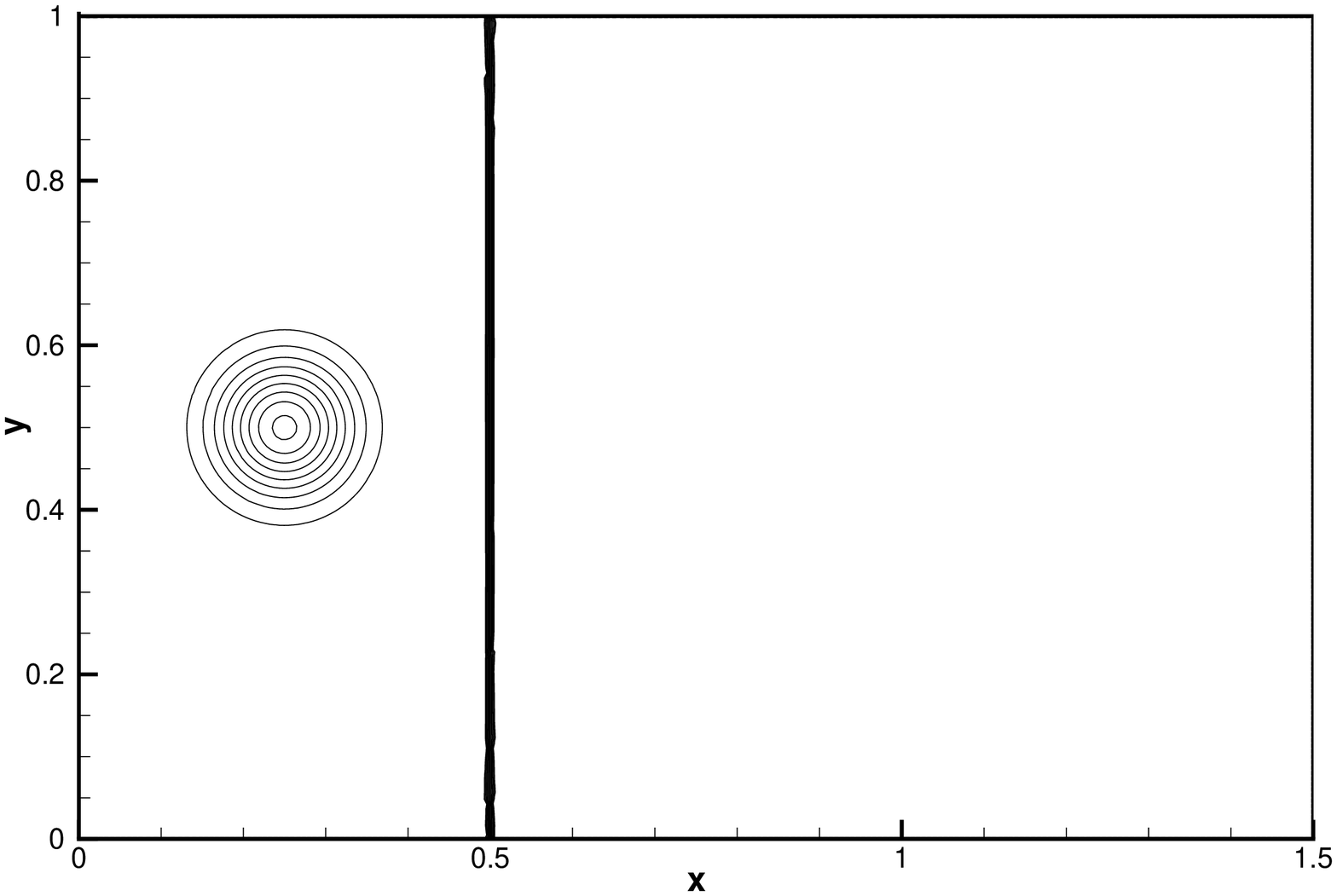}
\includegraphics[width=0.4\textwidth]{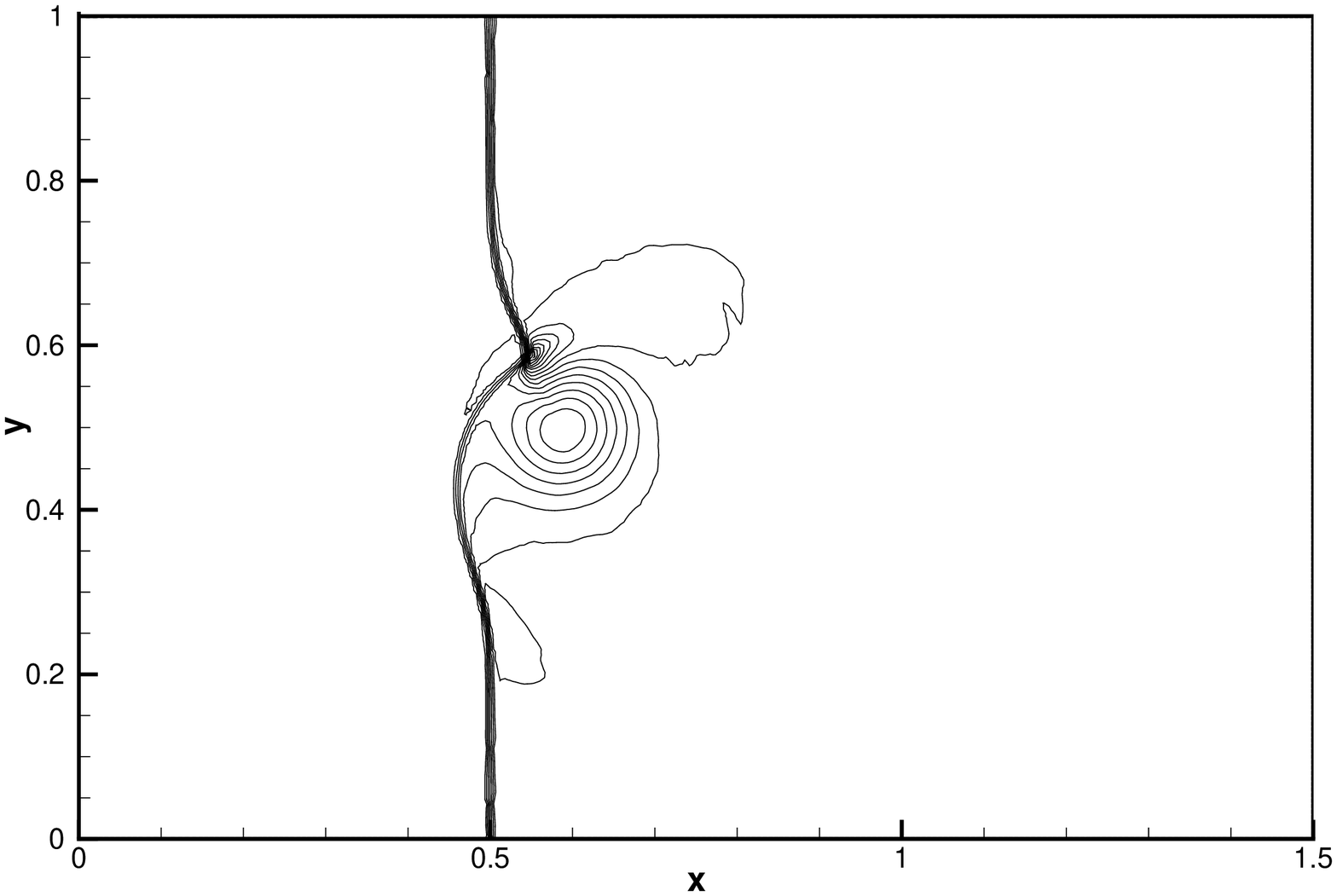}\\
\includegraphics[width=0.4\textwidth]{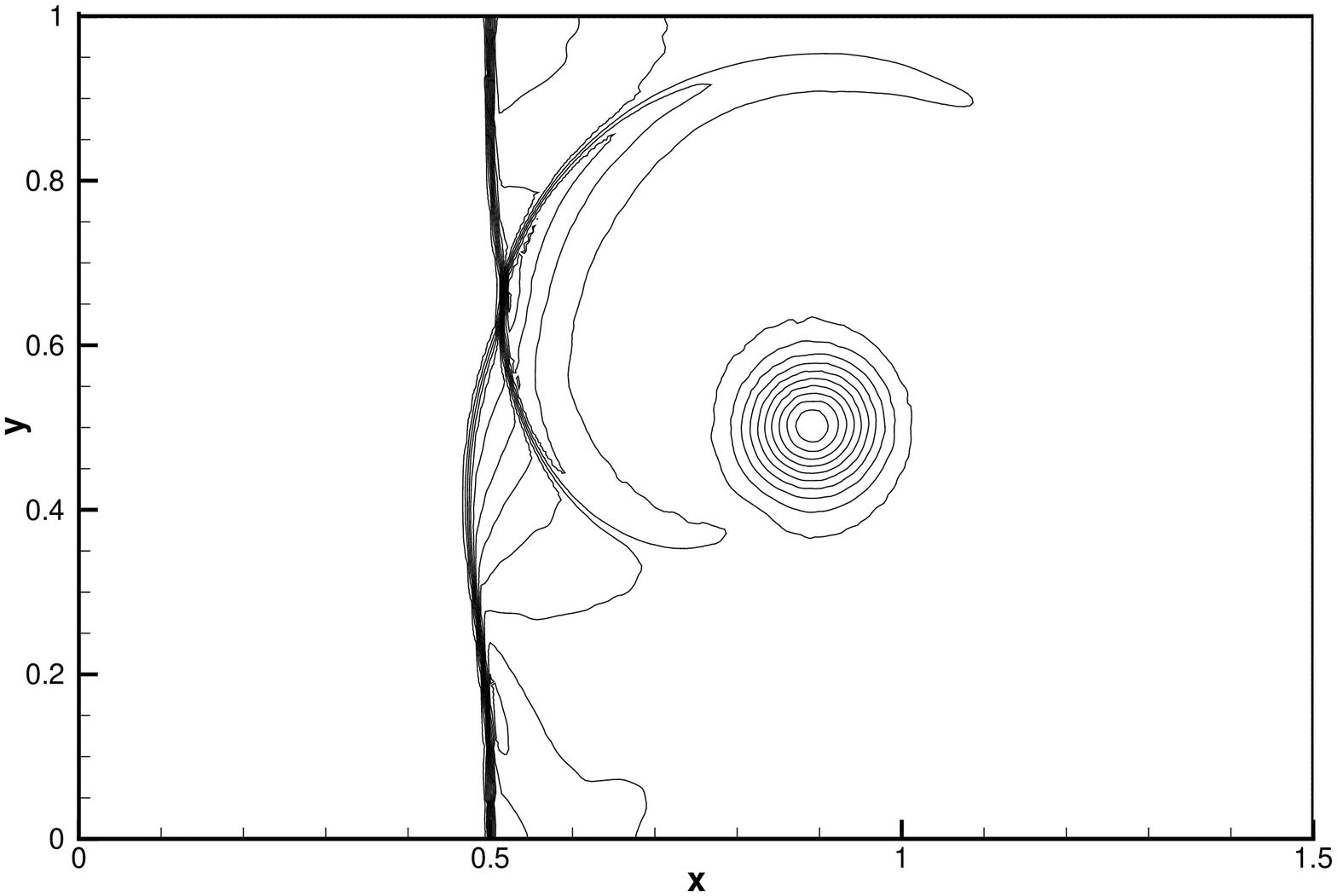}
\includegraphics[width=0.4\textwidth]{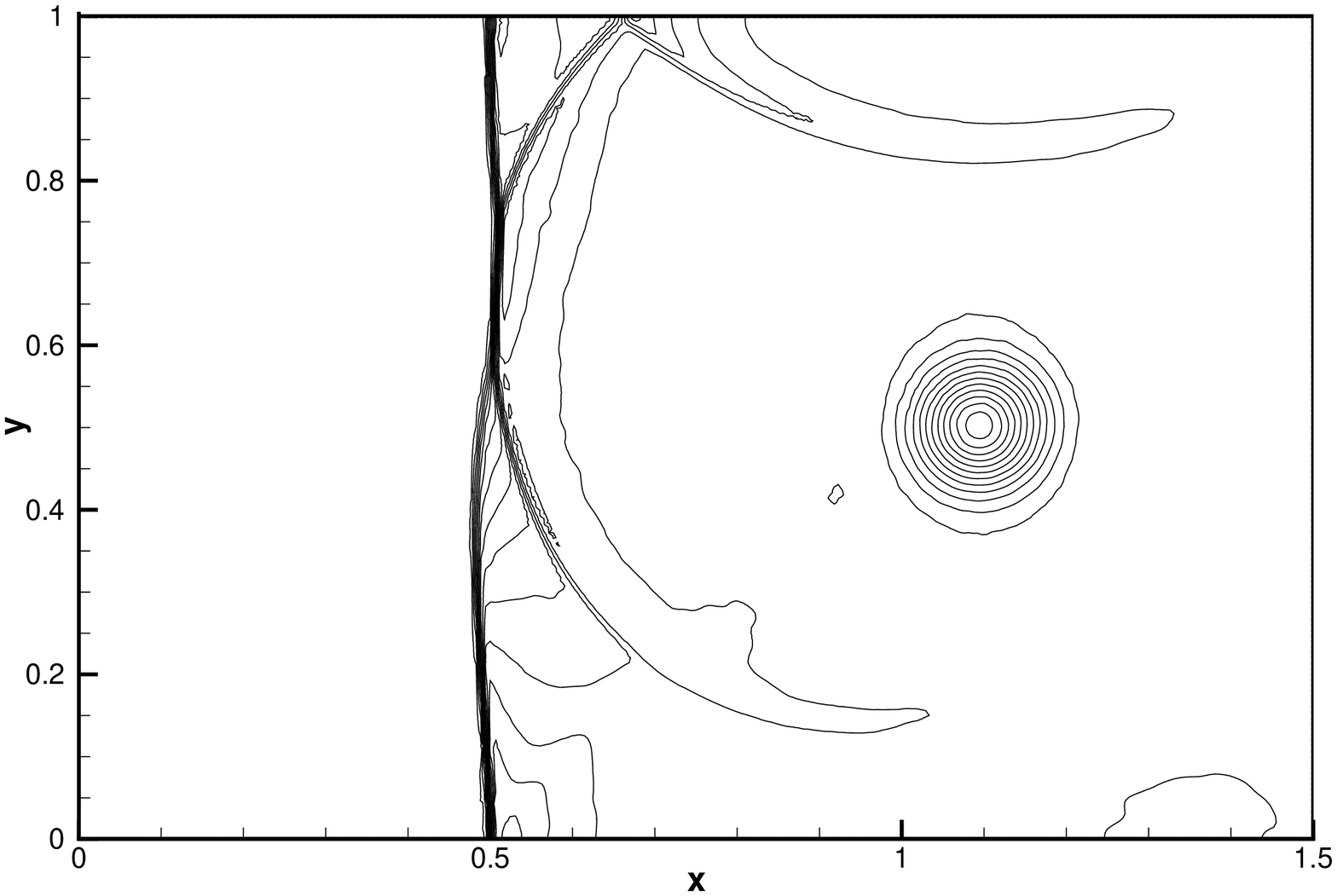}
\caption{\label{shock-vortex1}Shock vortex interaction: the pressure
distributions at $t=0.3$ and $0.8$ with mesh size $h=1/150$. }
\centering
\end{figure}

\begin{figure}[!h]
\centering
\includegraphics[width=0.375\textwidth]{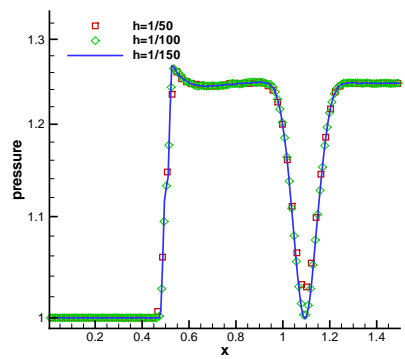}
\caption{\label{shock-vortex2}Shock vortex interaction: the pressure
distribution at $t=0.8$ along the horizontal symmetric line $y= 0.5$
with mesh size $h=1/50, 1/100$ and $1/150$. }
\end{figure}

The interaction between a stationary shock and a vortex for the
inviscid flow  is presented \cite{WENO2}. The computational domain is
taken to be $[0, 1.5]\times[0, 1]$. A stationary Mach $1.1$ shock is
positioned at $x=0.5$ and normal to the $x$-axis. The left upstream
state is $(\rho, u, v, p) = (Ma^2,\sqrt{\gamma}, 0, 1)$, where $Ma$
is the Mach number. A small vortex is obtained through a
perturbation on the mean flow with the velocity $(u, v)$,
temperature $T=p/\rho$ and entropy $S=\ln(p/\rho^\gamma)$, and the
perturbation is expressed as
\begin{align*}
&(\delta u,\delta v)=\kappa\eta e^{\mu(1-\eta^2)}(\sin\theta,-\cos\theta),\\
&\delta
T=-\frac{(\gamma-1)\kappa^2}{4\mu\gamma}e^{2\mu(1-\eta^2)},\delta
S=0,
\end{align*}
where $\eta=r/r_c$, $r=\sqrt{(x-x_c)^2+(y-y_c)^2}$, $(x_c,
y_c)=(0.25, 0.5)$ is the center of the vortex. Here $\kappa$ indicates
the strength of the vortex, $\mu$ controls the decay rate of the
vortex,  and $r_c$ is the critical radius for which the vortex has
the maximum strength. In the computation, $\kappa=0.3$, $\mu=0.204$,
and $r_c=0.05$. The reflected boundary conditions are used on the
top and bottom boundaries. The pressure distributions with mesh size
$h=1/150$ at $t=0, 0.3, 0.6$ and $0.8$ are shown in
Fig.\ref{shock-vortex1}. The detailed pressure distribution along
the center horizontal line with mesh size $h=1/50, 1/100$, and
$1/150$ at $t=0.8$ are shown in Fig.\ref{shock-vortex2}. This case
is tested without the detection of "trouble cell", which shows the
robustness of the weighted least square reconstruction to deal with
the flow with weak discontinuities.

\subsection{Double Mach reflection problem}
\begin{figure}[!h]
\centering
\includegraphics[width=0.449\textwidth]{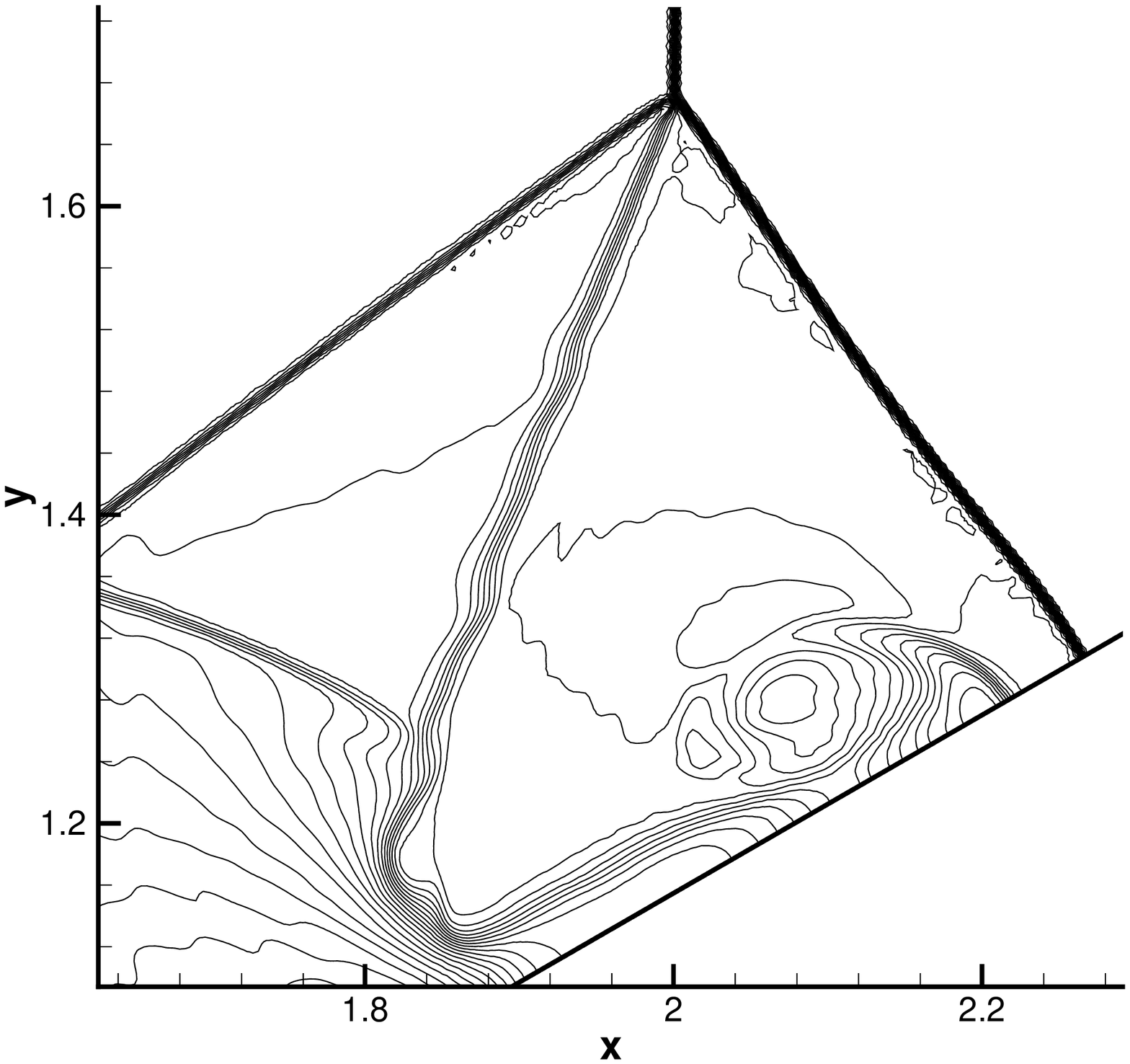}
\includegraphics[width=0.449\textwidth]{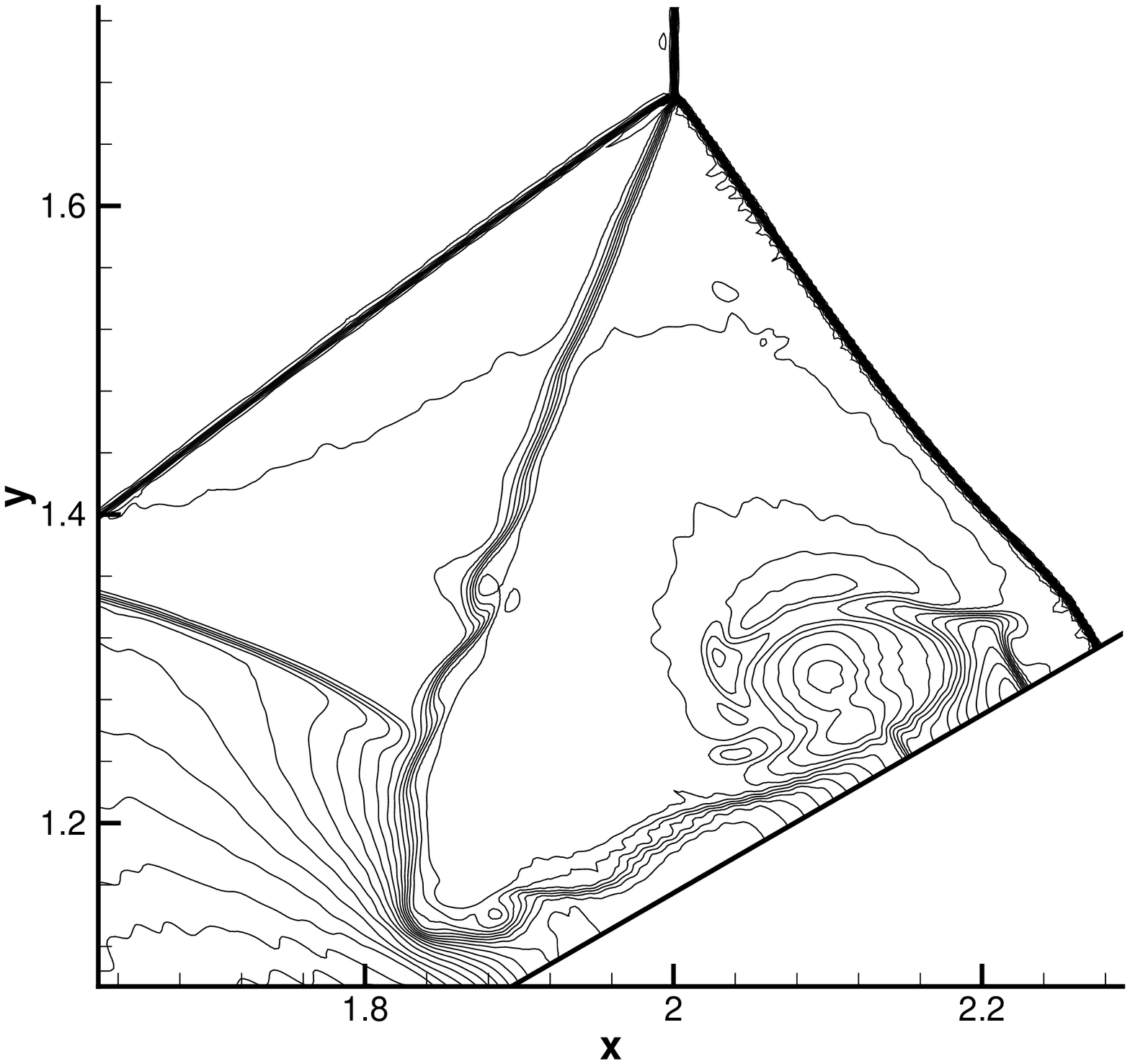}
\caption{\label{double-mach-2} Double Mach reflection: density
contours with the mesh size $1/240$, and $1/360$.}
\end{figure}

This problem was extensively studied by Woodward and Colella
 for the inviscid flow \cite{Case-Woodward}.
 A shock moves down in a tube which contains a $30^\circ$ wedge.
 The computational domain is shown in Fig.\ref{front-step2} with mesh
size $h=1/20$. The shock wave has a strength with Mach number $10$, which is
initially positioned at $x=0$. The initial pre-shock and post-shock
conditions are
\begin{align*}
(\rho, u, v, p)&=(8, 8.25, 0,
116.5),\\
(\rho, u, v, p)&=(1.4, 0, 0, 1).
\end{align*}
The reflective boundary conditions are used along the wedge, while
for the rest of bottom boundary, the exact post-shock condition is
imposed. At the top boundary, the flow variables are set to describe
the exact motion of the shock front along the wall. The inflow and
outflow boundary conditions are used at the entrance and the exit.
In this case, the weighted least square reconstruction is not
enough, and the shock detection technique is used to switch to the
second-order initial reconstruction. The density distributions with
mesh size $h=1/240$ and $1/360$ at $t=0.2$ are shown in
Fig.\ref{double-mach-2}. The compact scheme resolves the flow
structure under the triple Mach stem clearly.

\begin{figure}[!h]
\centering
\includegraphics[width=0.45\textwidth]{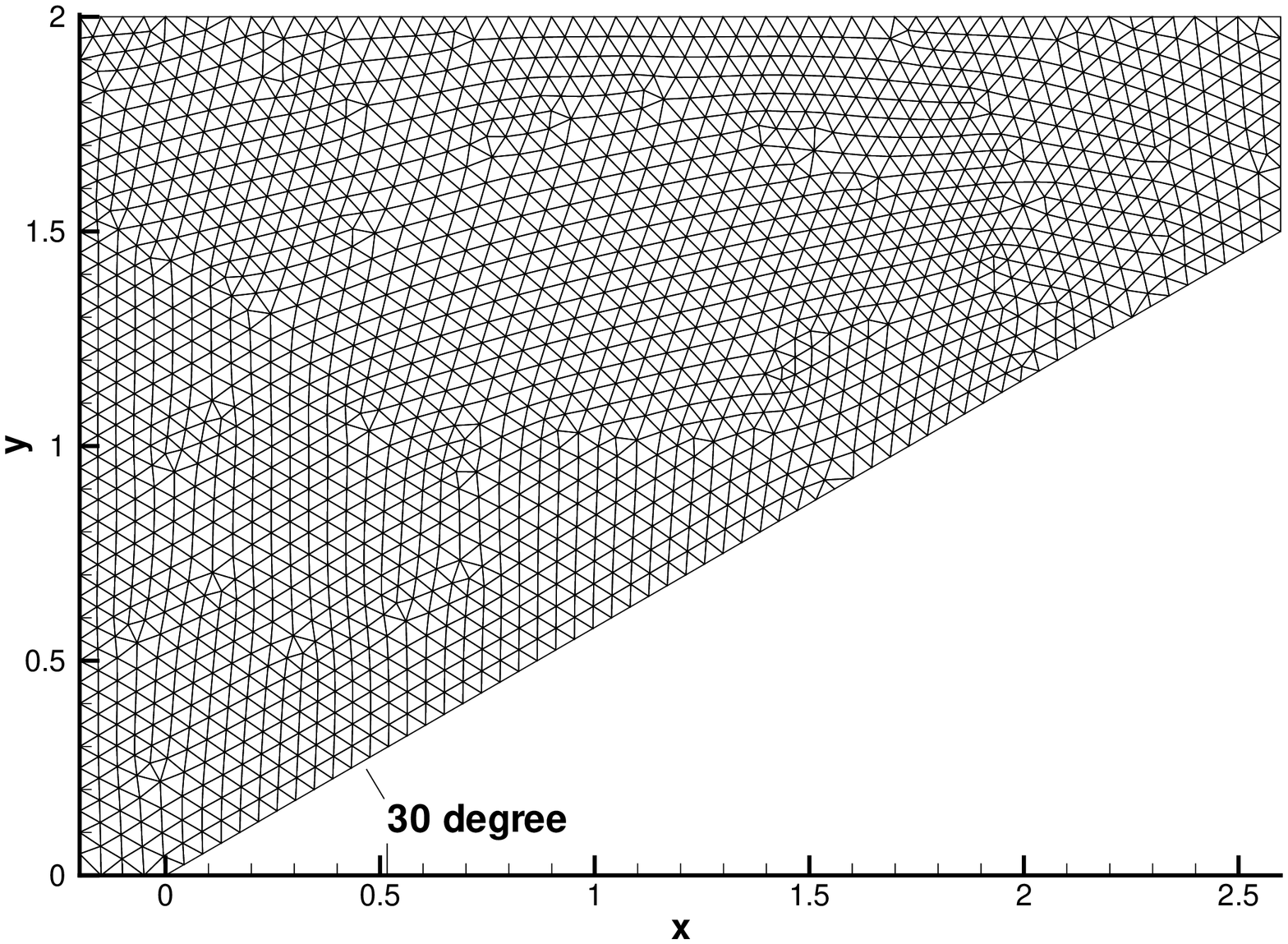}
\caption{\label{double-mach-1} Double Mach reflection: computational
domain with mesh size $h=1/20$.}
\end{figure}

\subsection{Mach step problem}
\begin{figure}[!h]
\centering
\includegraphics[width=0.25\textwidth]{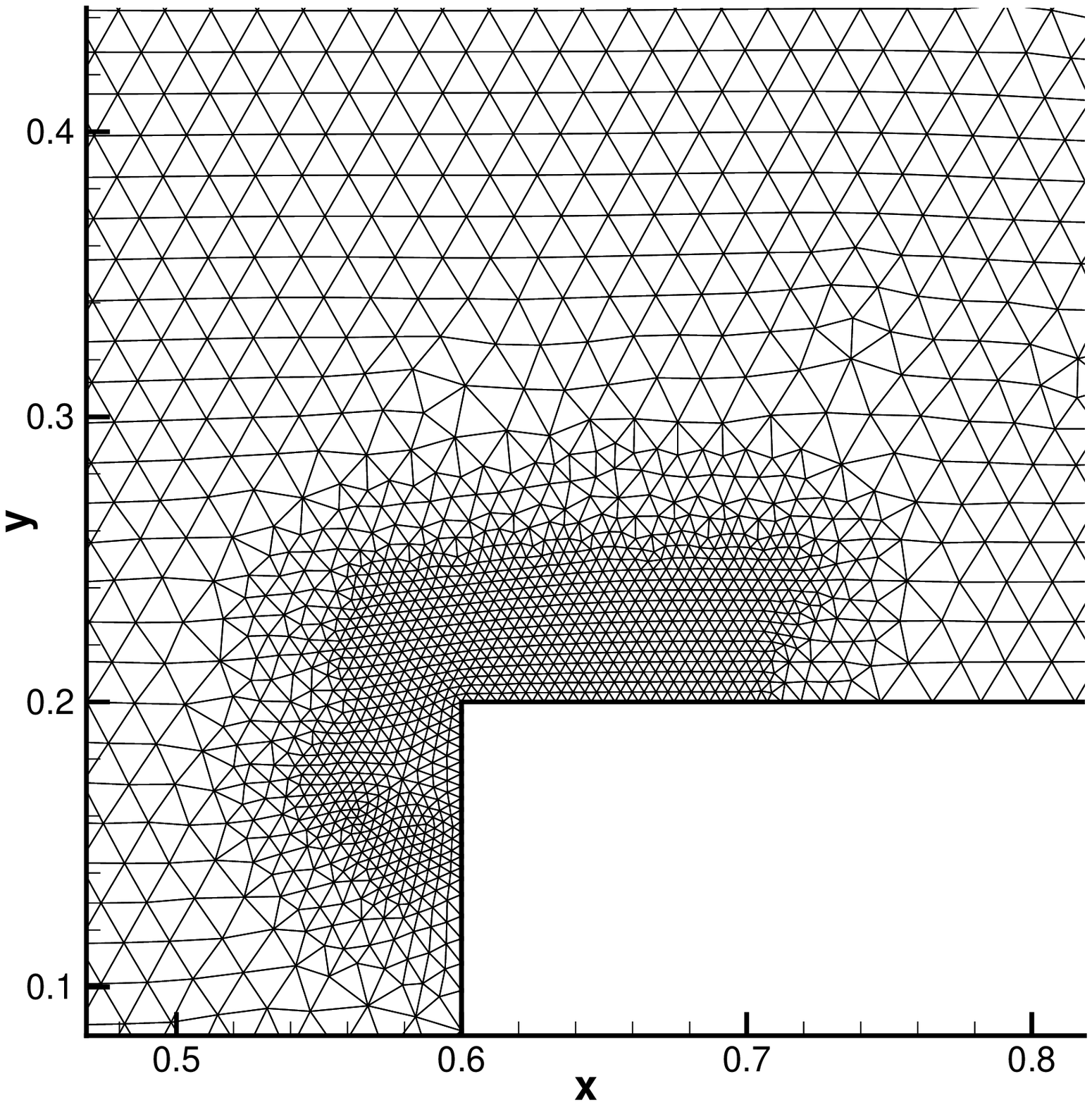}
\includegraphics[width=0.72\textwidth]{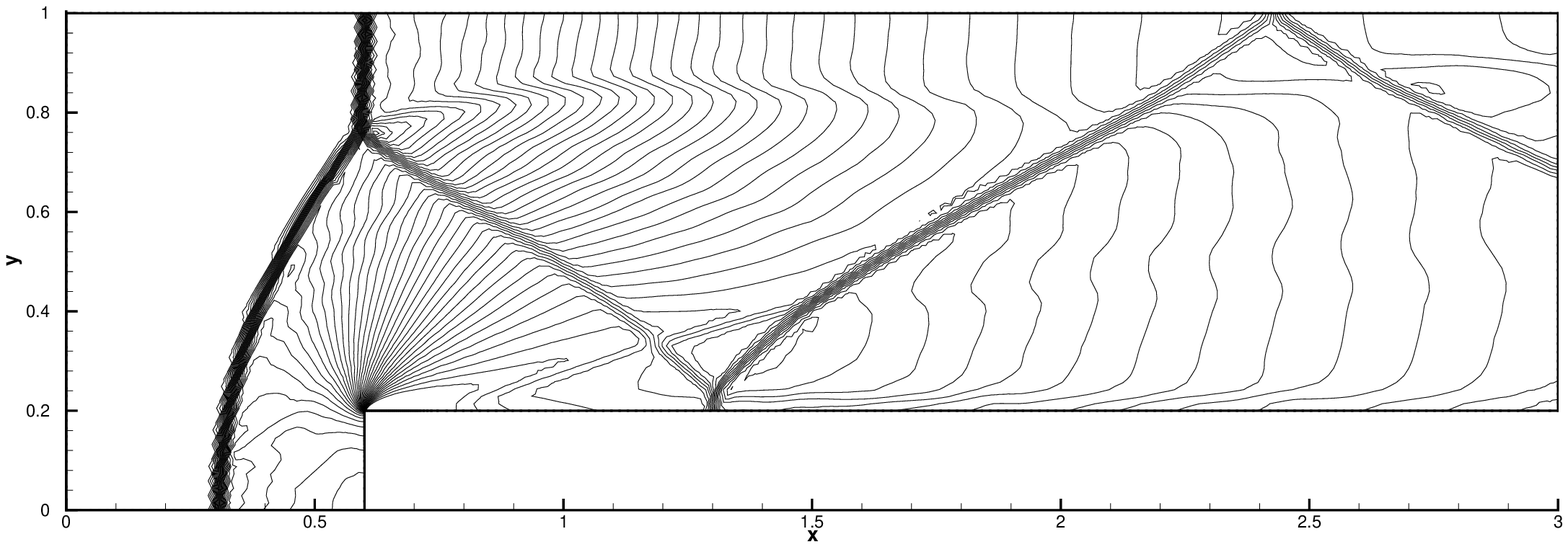}
\includegraphics[width=0.25\textwidth]{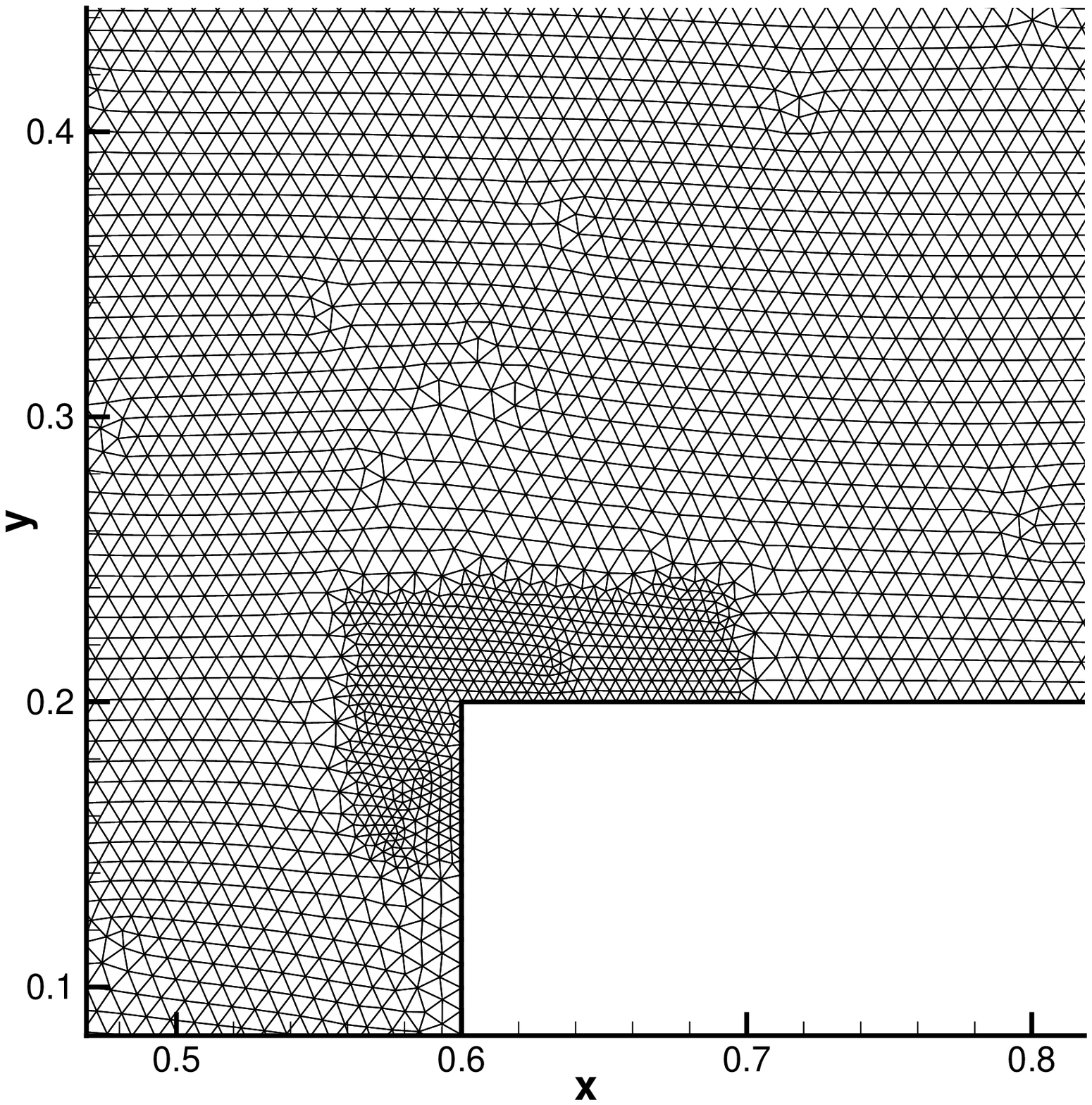}
\includegraphics[width=0.72\textwidth]{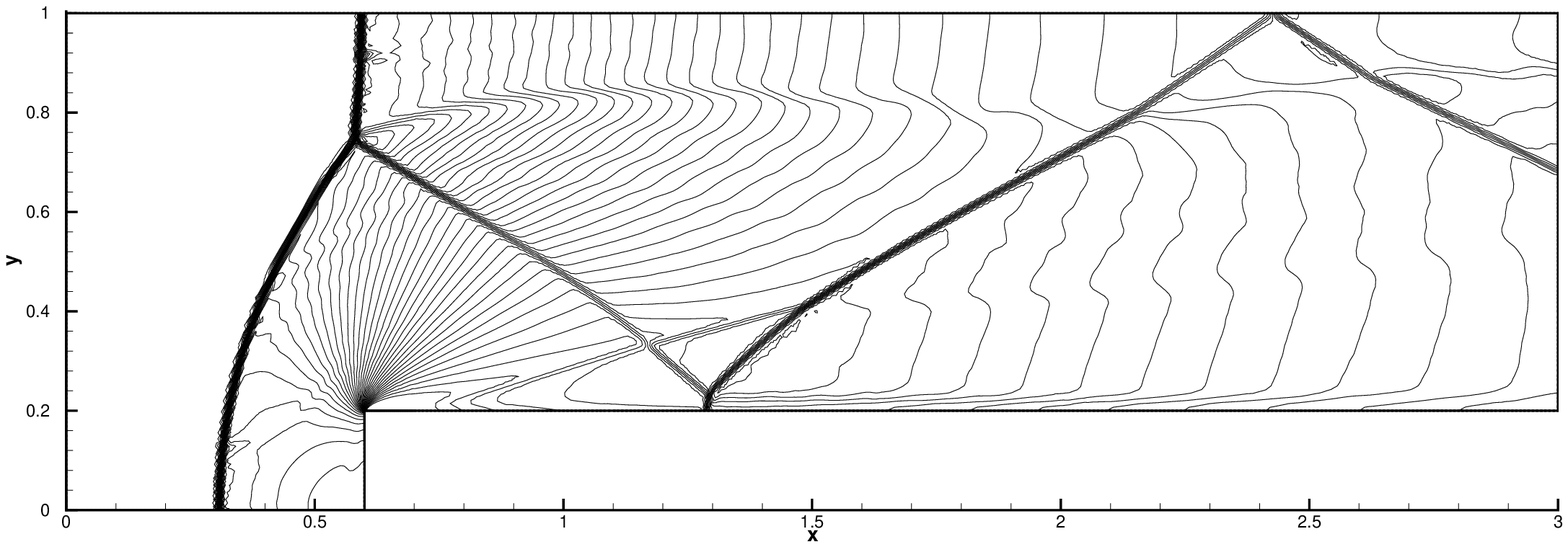}
\includegraphics[width=0.25\textwidth]{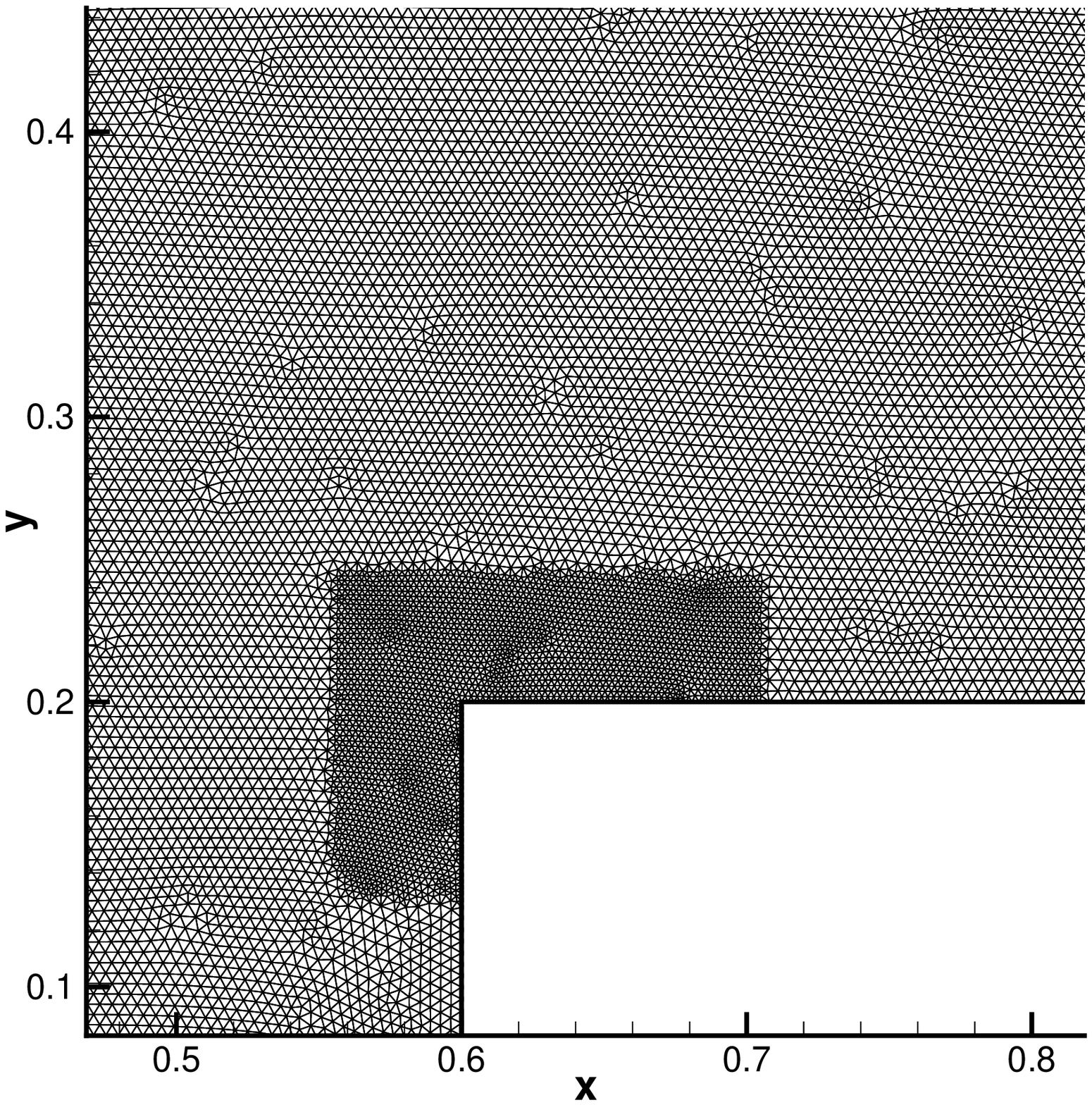}
\includegraphics[width=0.72\textwidth]{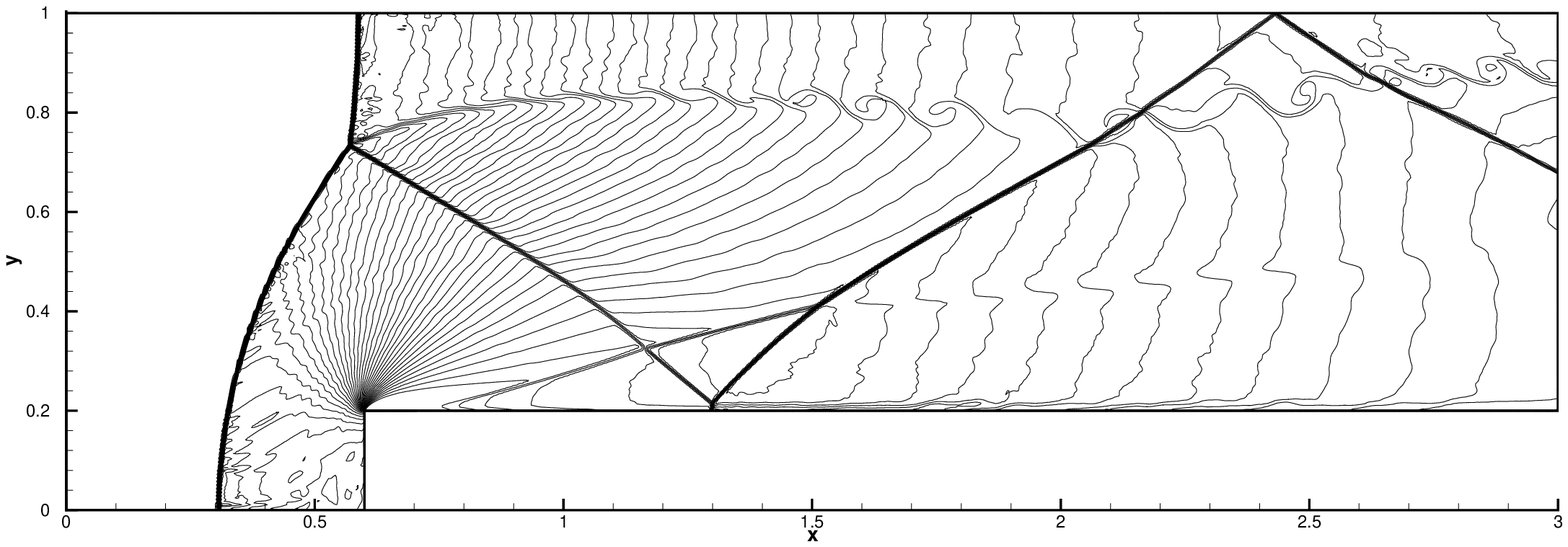}
\caption{\label{front-step2}Mach step problem: density contours
with the mesh size $h=1/60, 1/120$, and $1/240$.}
\end{figure}

The Mach step problem was again studied extensively by Woodward and
Colella \cite{Case-Woodward} for the inviscid flow. The
computational domain $[0,3]\times[0,1]$ is shown in
Fig.\ref{front-step2}, which is covered by unstructured mesh with mesh size $h=1/20$.
The Mach step is
located at $x=0.6$ with height $0.2$ in the tunnel. Initially, a
right-moving flow with Mach $3$ is imposed in the whole computational domain.
The reflective boundary conditions
are used along the walls of the tunnel, and inflow and outflow
boundary conditions are used at the entrance and the exit.
The corner of the step is the center of a rarefaction fan, which is  a
singularity point. To minimize the
numerical error generated at the corner, the meshes near
the corner are refined, shown in Fig.\ref{front-step2}.
In this case, the weighted least square reconstruction is not enough,
and the shock detection technique is used again to switch to the second-order reconstruction.
The density distributions
with $h=1/60, 1/120$, and $1/240$ at $t=4$ are presented in
Fig.\ref{front-step2}. With the mesh refinement, the resolution is
improved, especially for the slip line started from the triple point.

\begin{figure}[!h]
\centering
\includegraphics[width=0.6\textwidth]{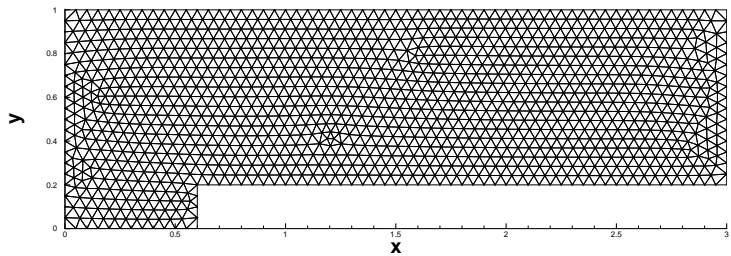}
\caption{\label{front-step1}Mach step problem: the computational
domain with mesh size $h=1/20$.}
\end{figure}

\subsection{Lid-driven cavity flow}
\begin{figure}[!h]
\centering
\includegraphics[width=0.44\textwidth]{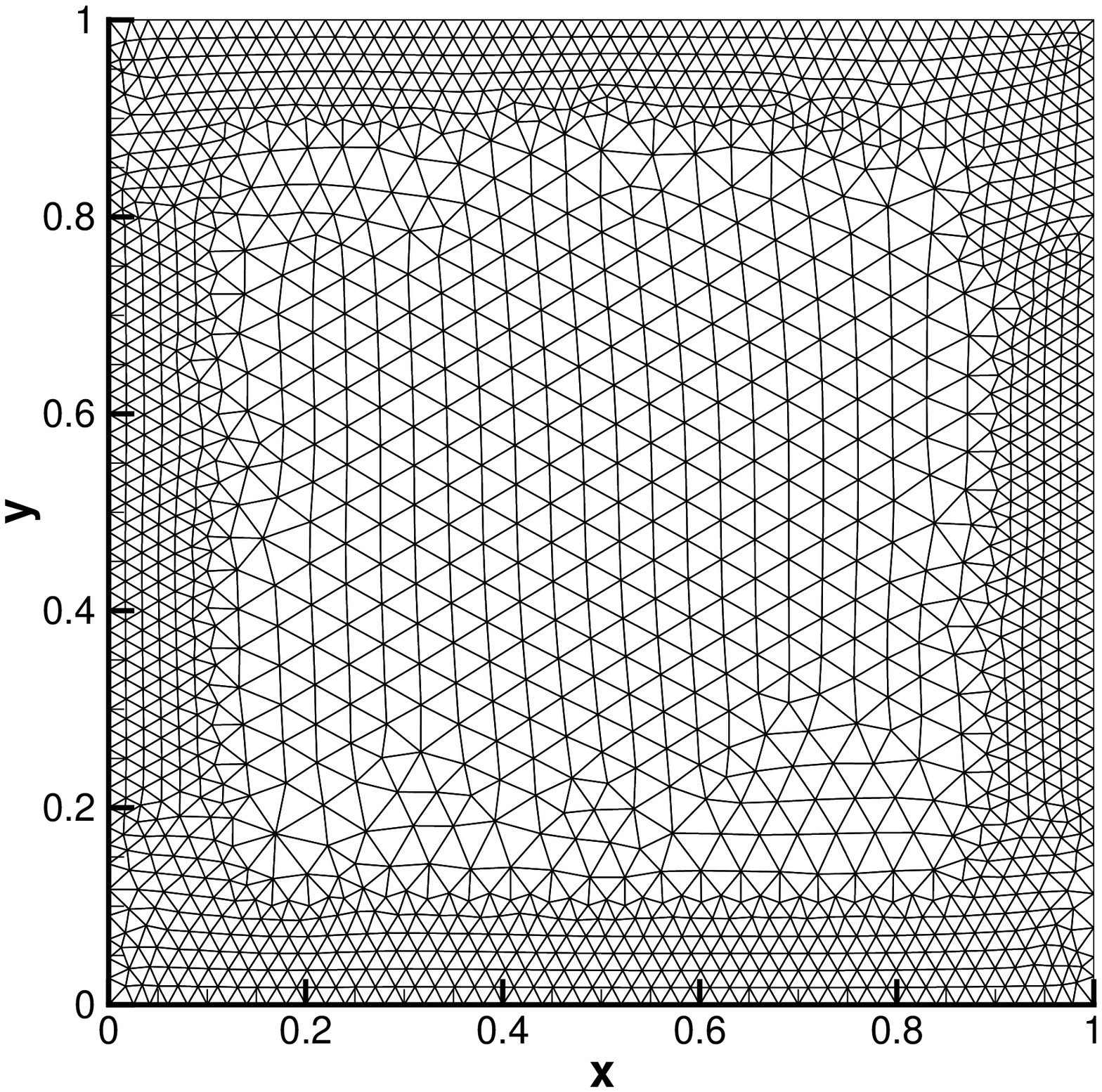}
\includegraphics[width=0.44\textwidth]{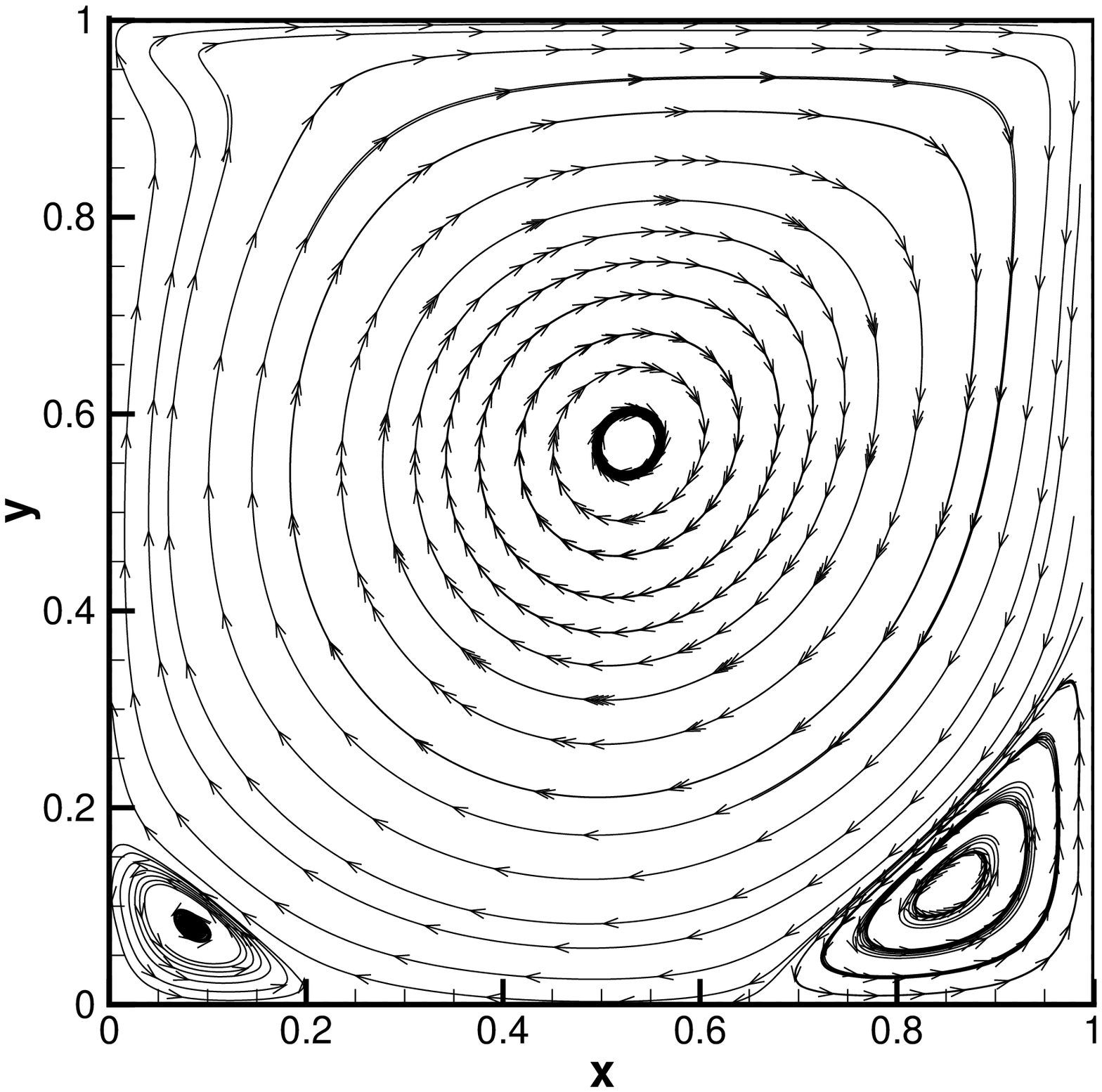}
\caption{\label{cavity-1} Lid-driven cavity flow: mesh and
streamlines for the compact gas-kinetic scheme with $Re=1000$.}
\end{figure}

\begin{figure}[!h]
\centering
\includegraphics[width=0.424\textwidth]{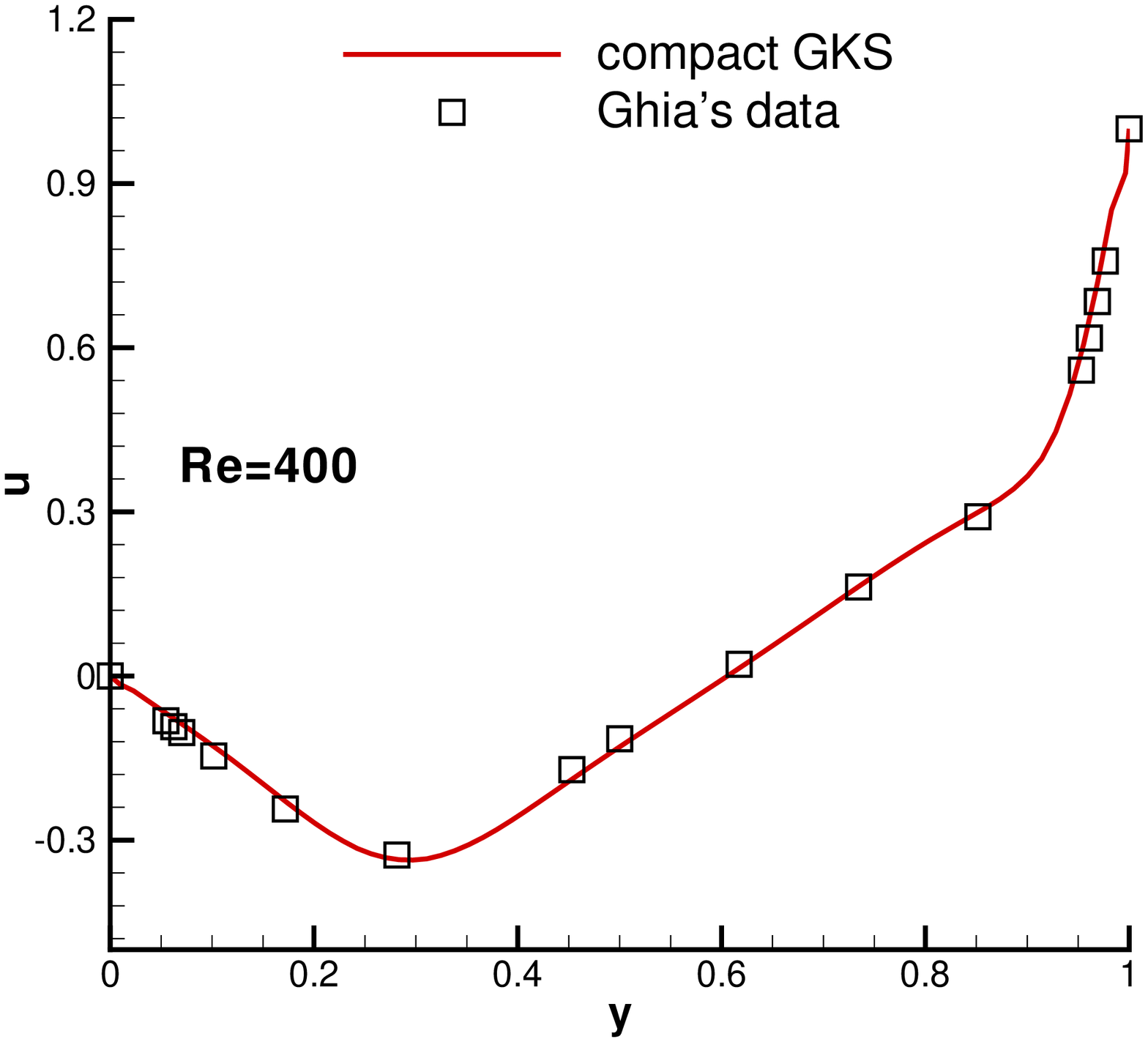}
\includegraphics[width=0.424\textwidth]{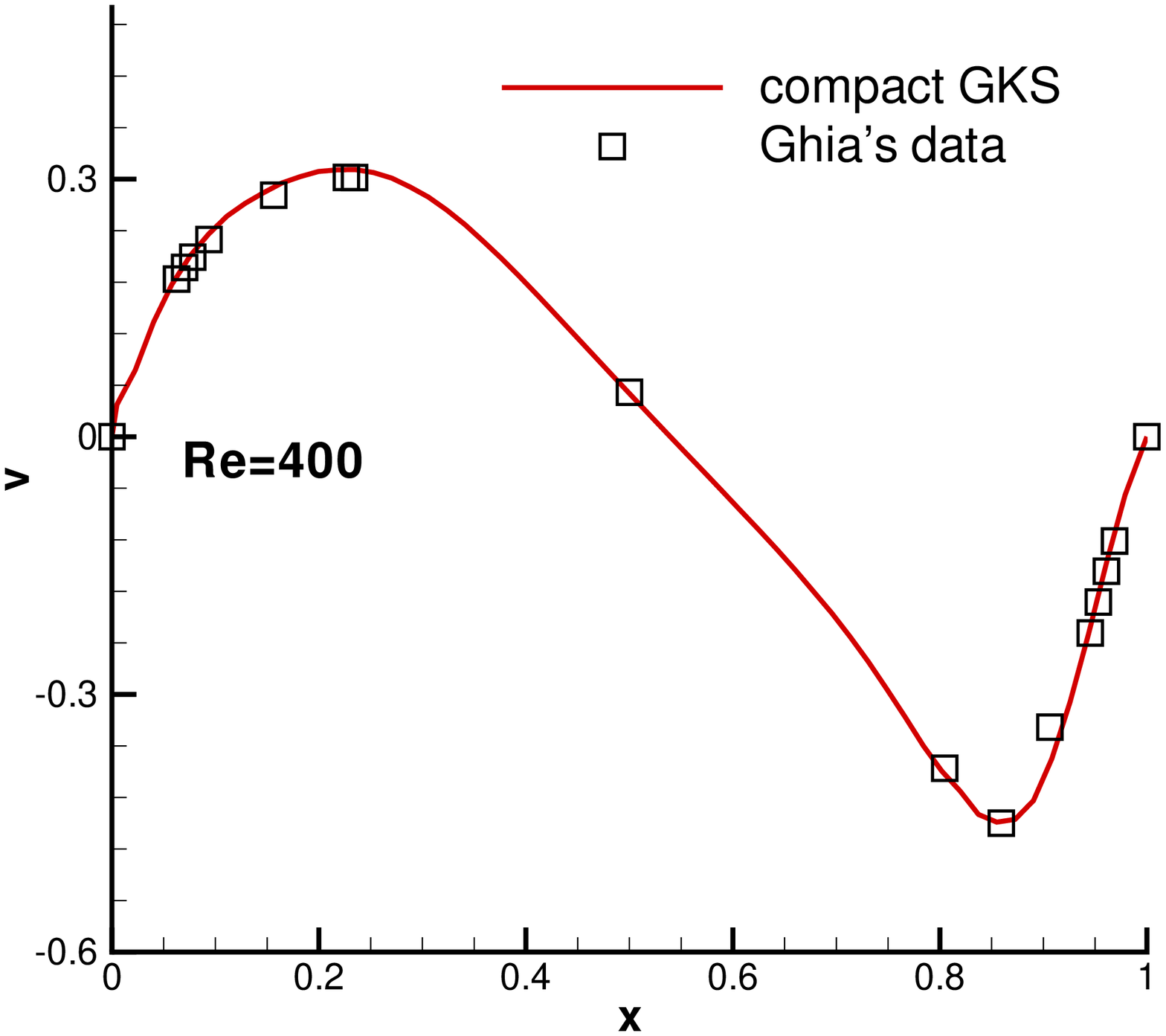}\\
\includegraphics[width=0.424\textwidth]{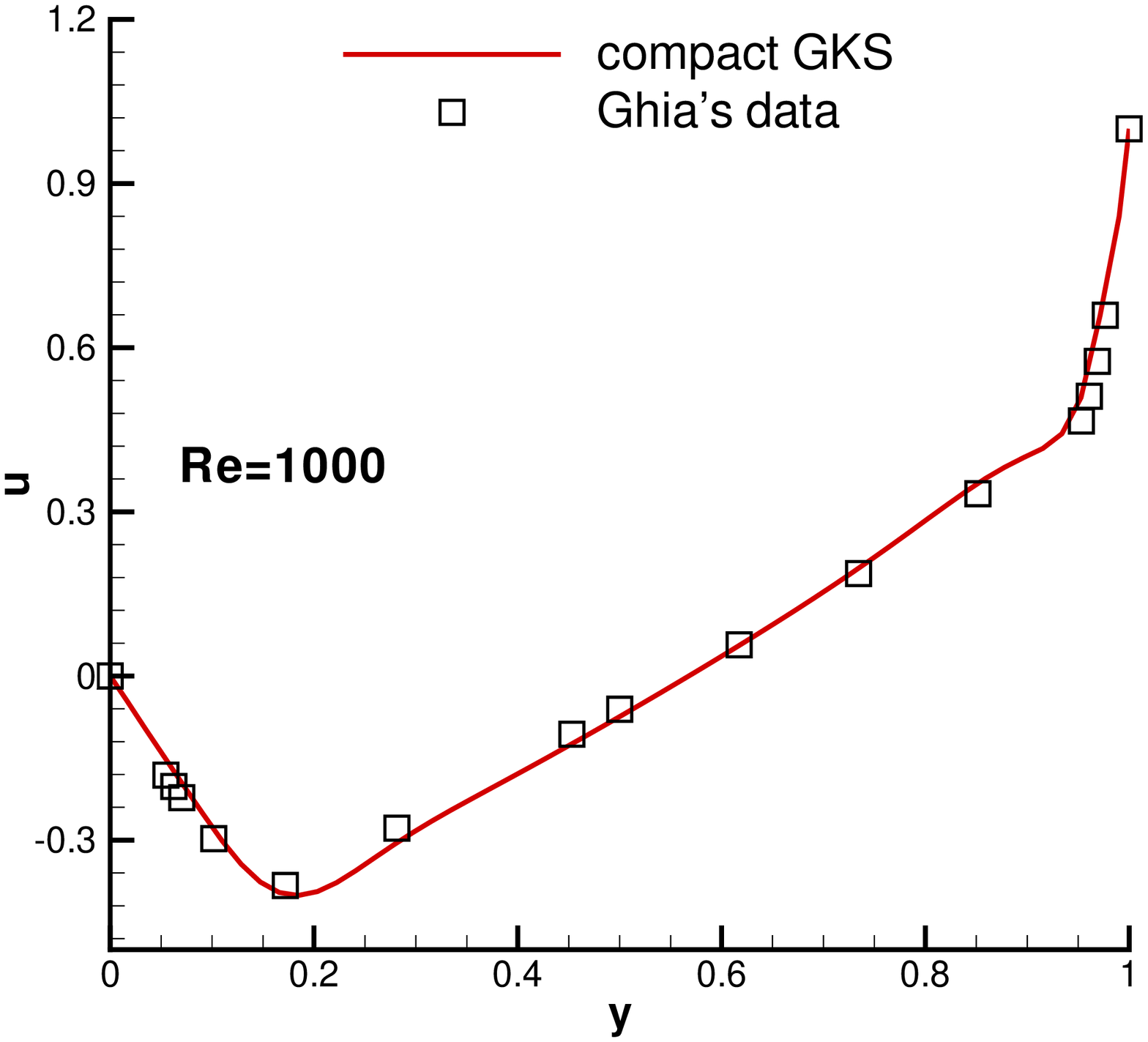}
\includegraphics[width=0.424\textwidth]{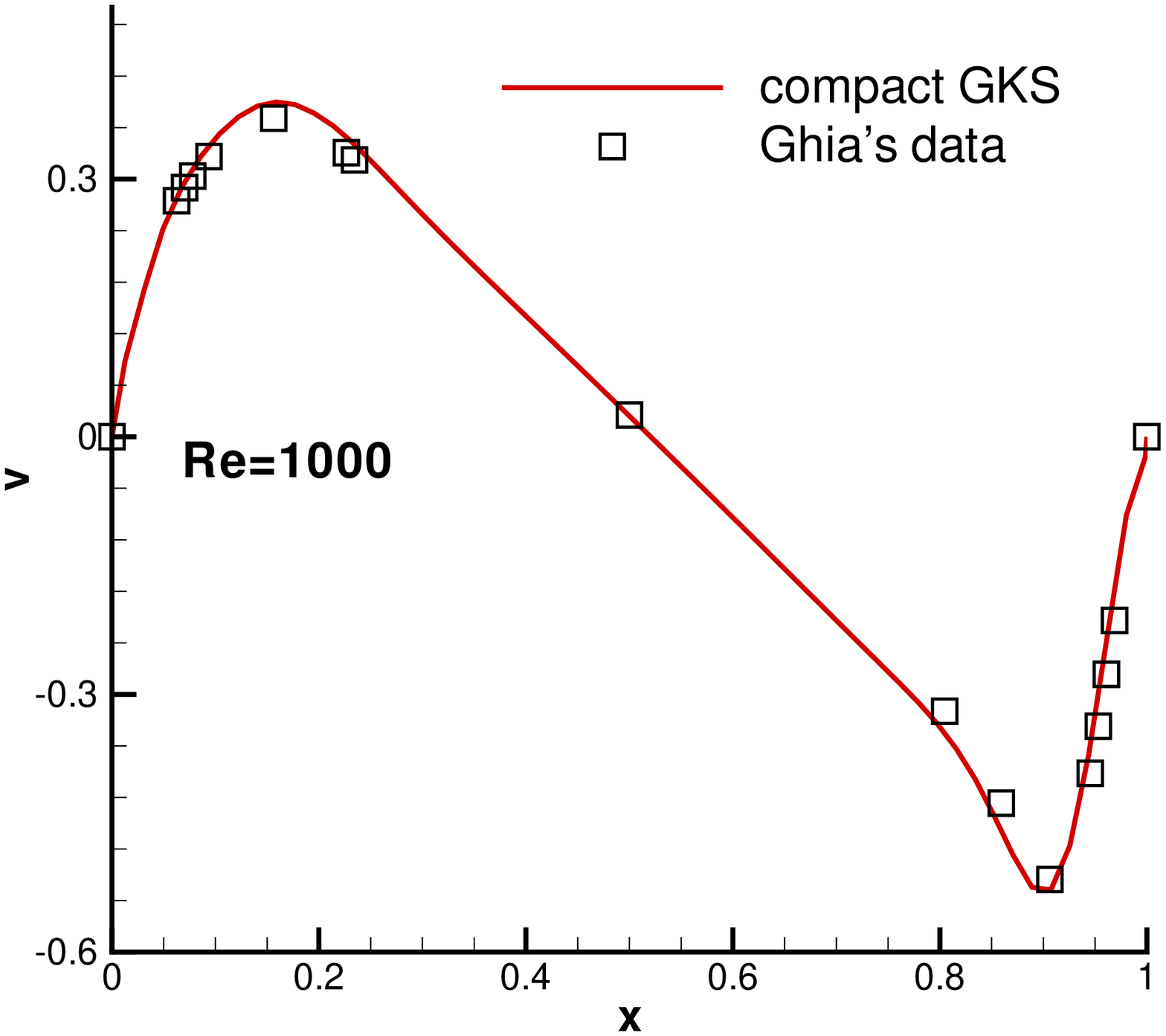}\\
\includegraphics[width=0.424\textwidth]{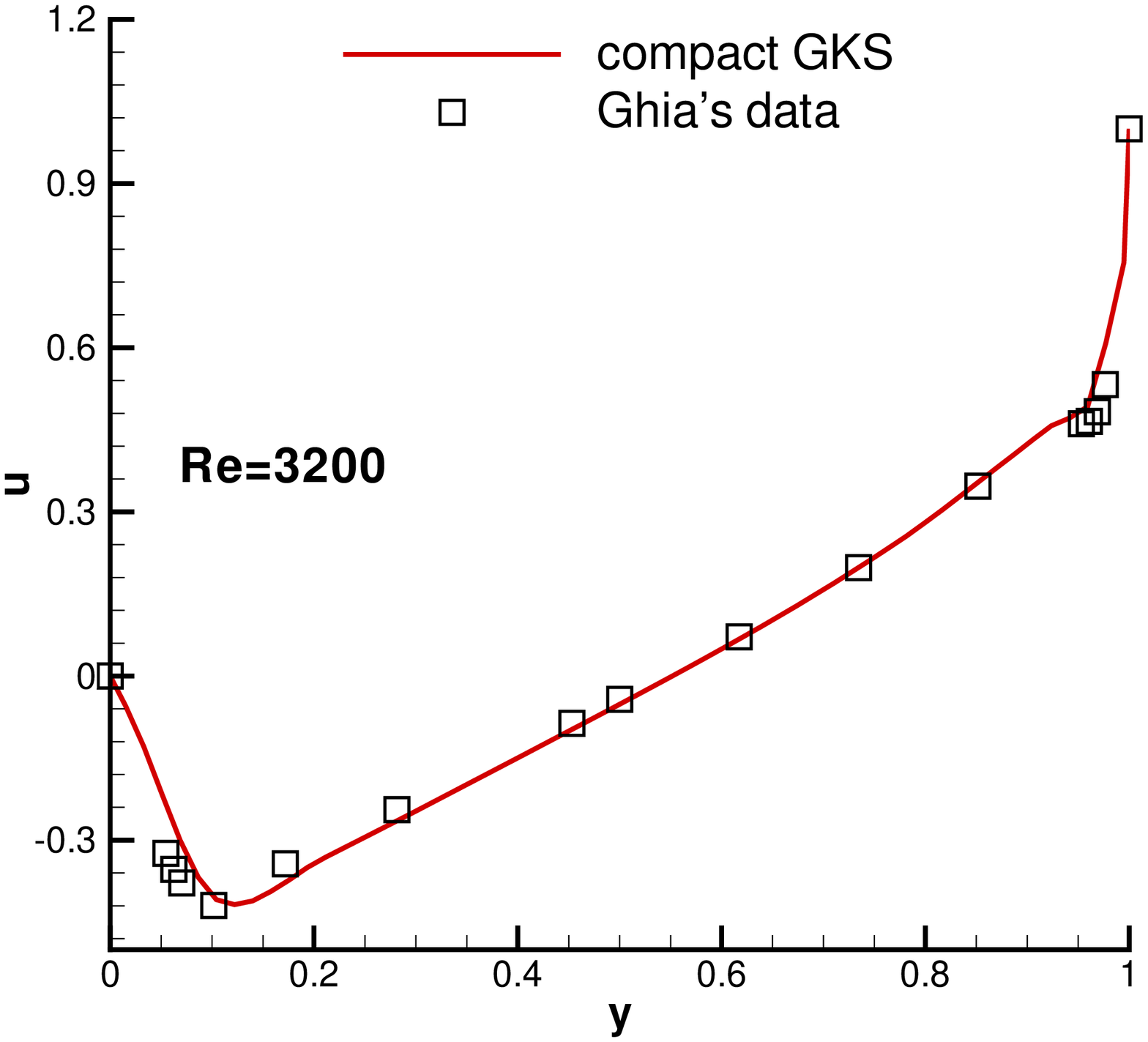}
\includegraphics[width=0.424\textwidth]{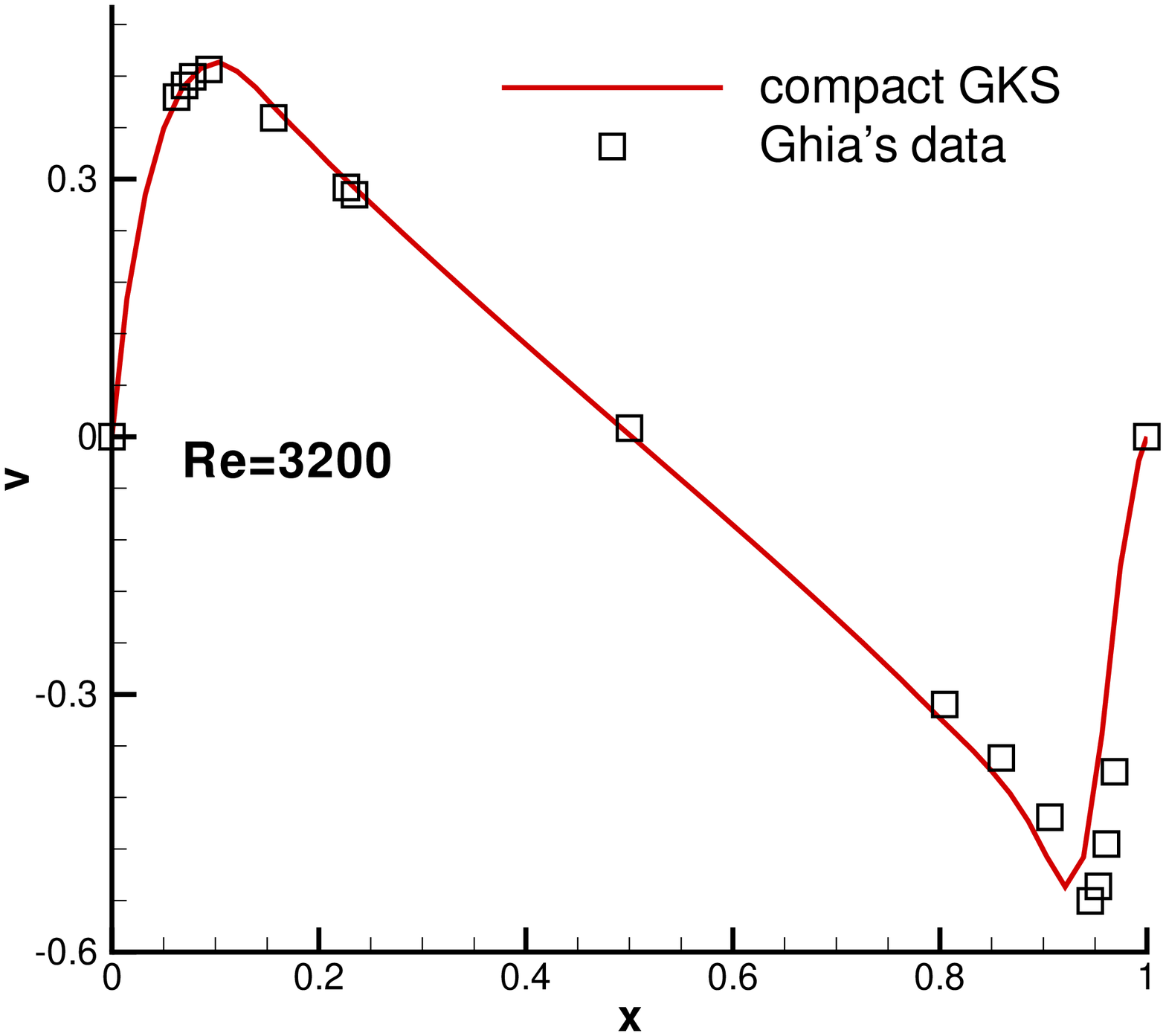}
\caption{\label{cavity-2} Lid-driven cavity flow:  $U$-velocities
along vertical centerline line and $V$-velocities along horizontal
centerline with $Re=400, 1000$ and $3200$. The reference data is
from Ghia \cite{Case-Ghia}.}
\end{figure}

The lid-driven cavity problem is one of the most important
benchmarks for validating incompressible or low speed Navier-Stokes
flow solvers. The fluid is bounded by a unit square and driven by a
uniform translation of the top boundary. In this case,  the gas has a  specific heat ratio
$\gamma=5/3$ and the up wall is
moving with a speed of Mach number $Ma=0.15$.
Isothermal and nonslip boundary conditions are imposed.
The computational domain $[0, 1]\times[0, 1]$ with unstructured mesh is presented
in Fig.\ref{cavity-1}, where mesh size  are
$h=1/25$ for the inner cells and  $h=1/50$ near the walls.  Numerical
simulations are conducted for three Reynolds numbers $Re=400, 1000$
and $3200$. The streamlines with $Re=1000$ for the compact scheme
are shown in Fig.\ref{cavity-1}. The results of
$U$-velocities along the center vertical line, $V$-velocities along
the center horizontal line, and the benchmark data \cite{Case-Ghia}
are shown in Fig.\ref{cavity-2} for $Re= 400, 1000, 3200$. The
simulation results match well with the benchmark data.

\begin{figure}[!h]\centering
\includegraphics[width=0.88\textwidth]{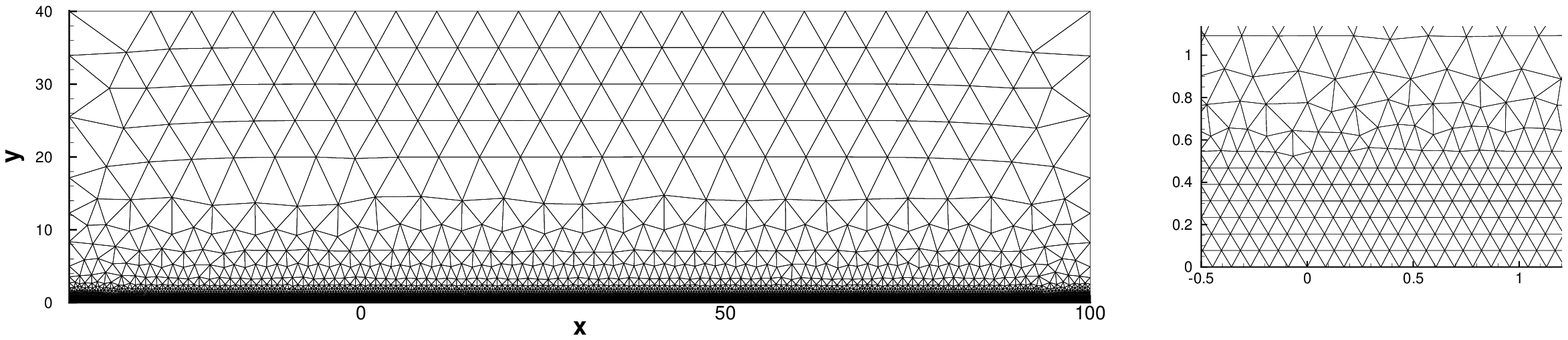}
\includegraphics[width=0.88\textwidth]{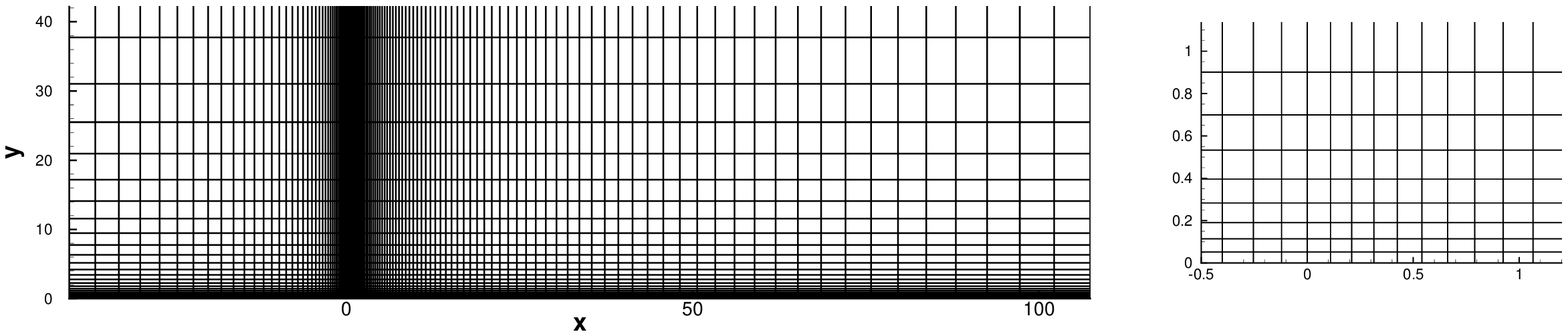}
\caption{\label{boundary-layer1}Laminar boundary layer computation: the triangular
and rectangular meshes.}
\end{figure}

\begin{figure}[!h]\centering
\includegraphics[width=0.47\textwidth]{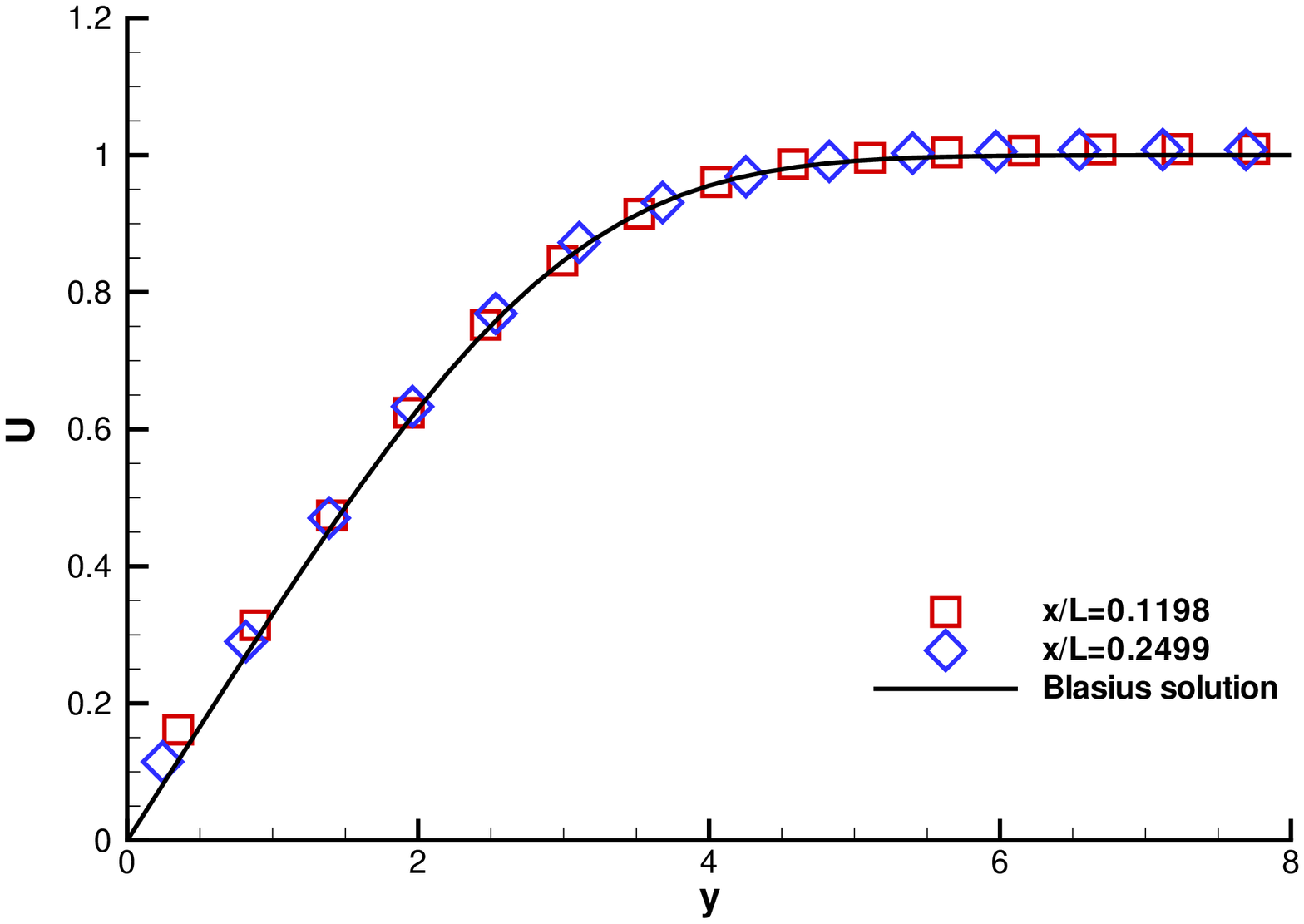}
\includegraphics[width=0.47\textwidth]{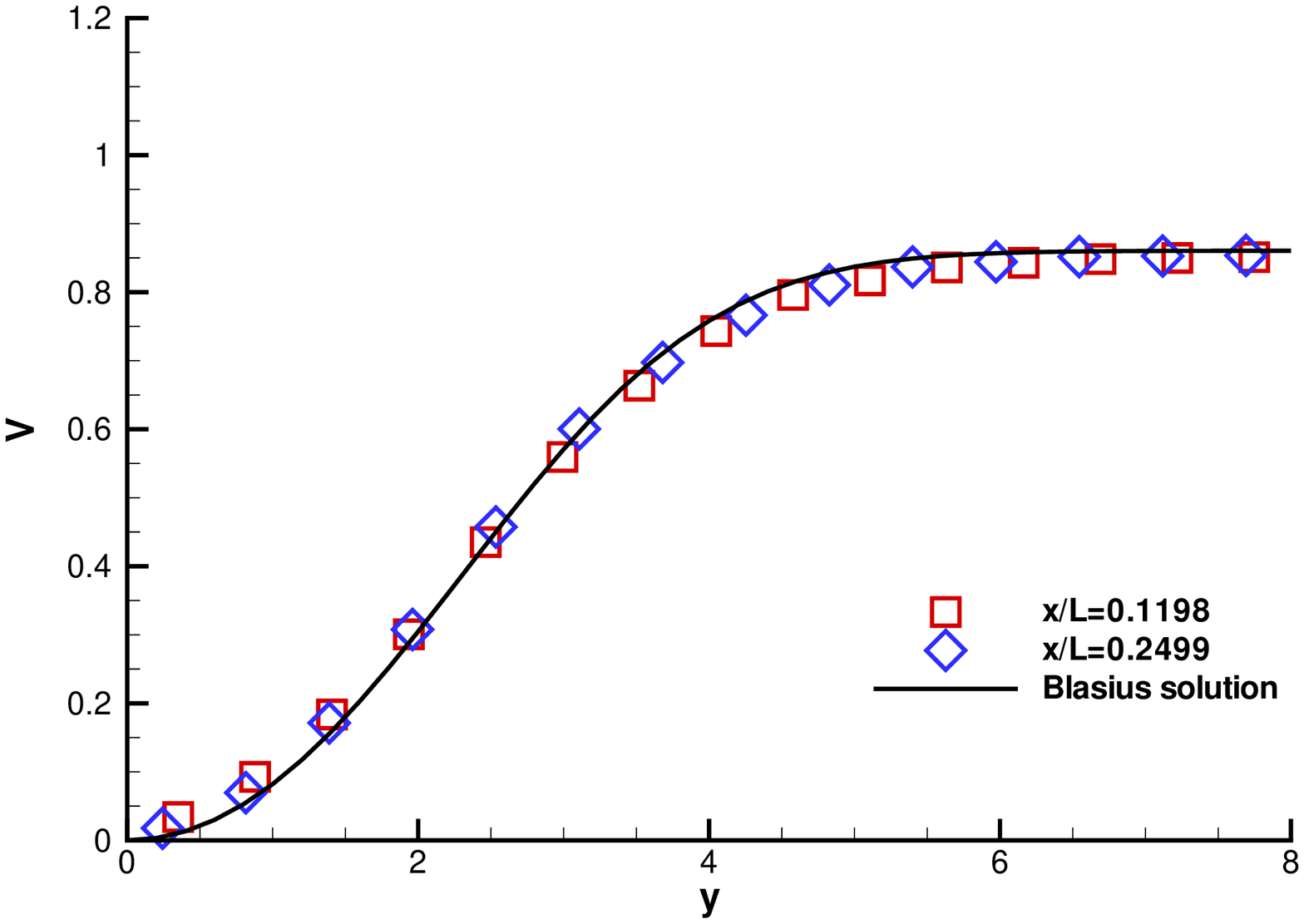}
\caption{\label{boundary-layer2}Laminar boundary layer solution from compact scheme with triangular mesh: the
non-dimensional velocity $U$ and $V$. }
\includegraphics[width=0.47\textwidth]{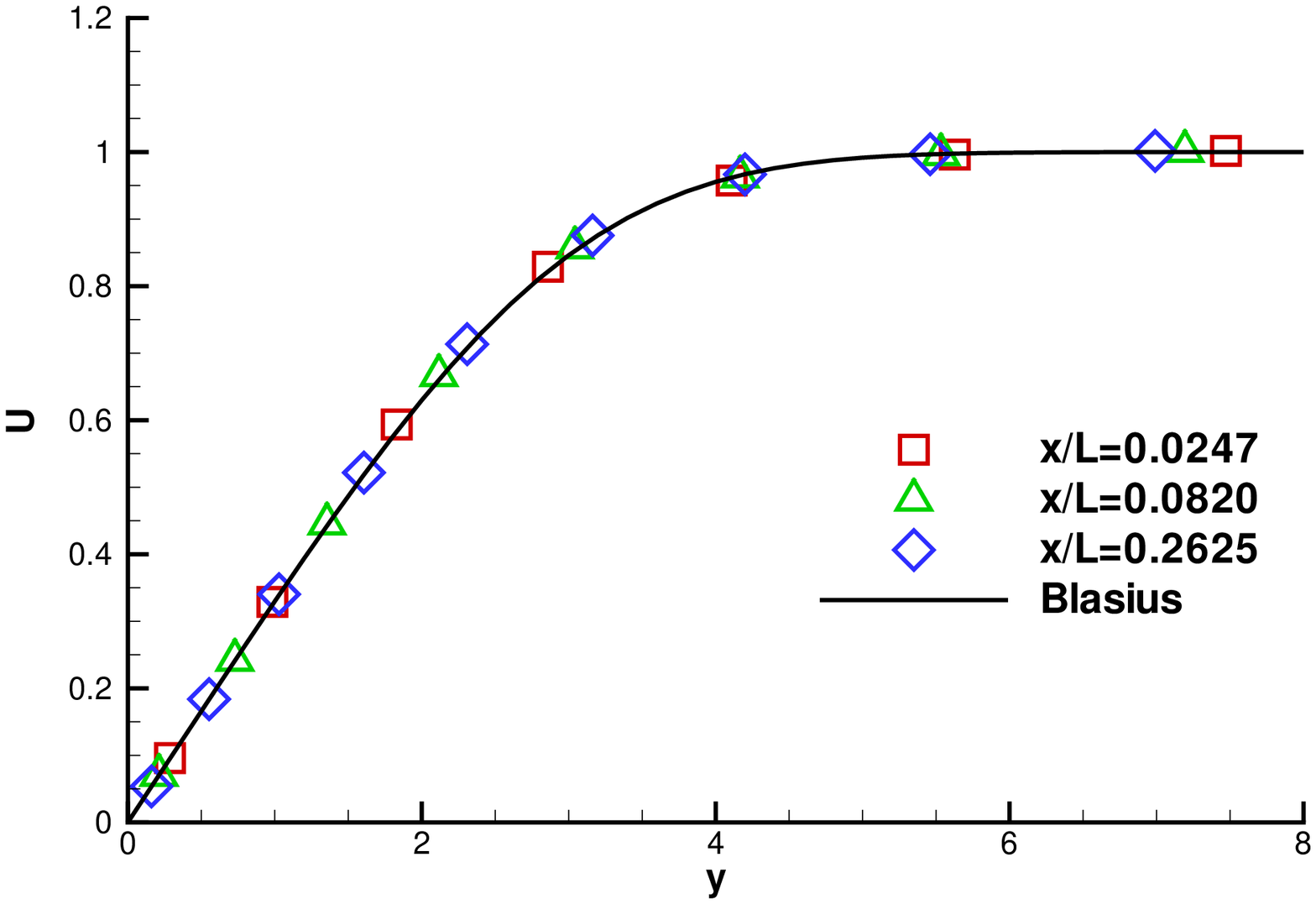}
\includegraphics[width=0.47\textwidth]{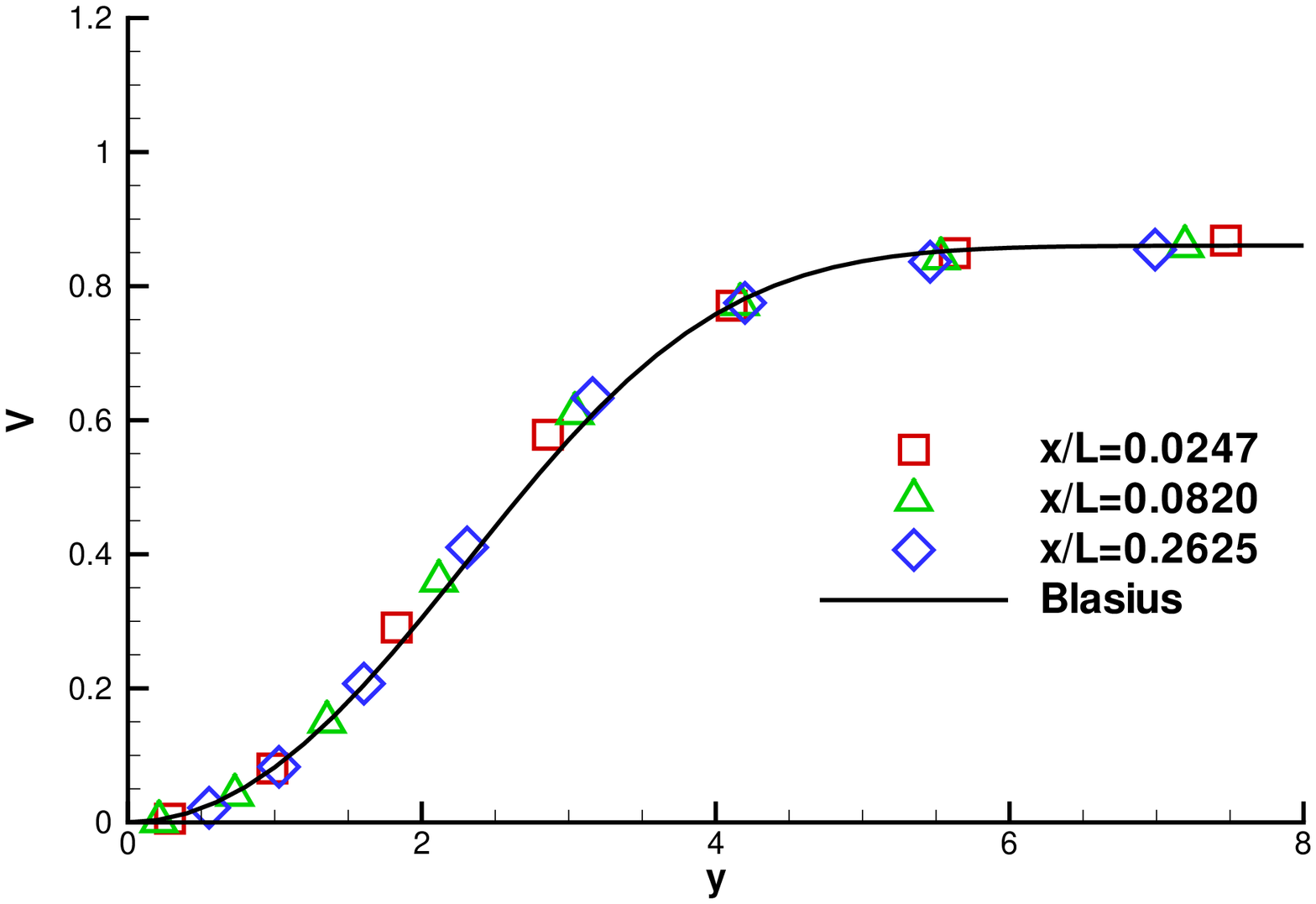}
\caption{\label{boundary-layer3}Laminar boundary layer solution from compact scheme with rectangular mesh: the
non-dimensional velocity $U$ and $V$. }
\end{figure}

\subsection{Laminar boundary layer}
A laminar boundary layer is tested over a flat plate with length
$L=100$. The Mach number of the free-stream is $Ma=0.15$ and the
Reynolds number is $Re=U_{\infty}L/\nu=10^5$, $\nu$ is the viscous
coefficient. This case is tested with the compact scheme for the
both triangle mesh and rectangular mesh.
 Fig.\ref{boundary-layer1} presents both triangular  and
rectangular meshes, with an enlarged view of meshes near the boundary.
The non-slip adiabatic boundary
condition at the plate is used and a symmetry condition is imposed
at the bottom boundary before the flat plate. The non-reflecting
boundary condition based on the Riemann invariants is adopted for the
other boundaries. The non-dimensional velocity $U$ and $V$ at different locations are given
in Fig.\ref{boundary-layer2} for the
triangular mesh and Fig.\ref{boundary-layer3} for the rectangular
mesh. In all locations, the numerical solutions match with  the exact Blasius
solutions very well. Here the boundary layer can be resolved by six or seven
mesh points. The solutions show the good performance of the compact
scheme for the Navier-Stokes solutions with unstructured mesh.

\subsection{Viscous shock tube problem}
This problem was introduced in \cite{Case-Daru} to test the
performances of different schemes for viscous flows. In this case,
an ideal gas is at rest in a two-dimensional unit box
$[0,1]\times[0,1]$. A membrane located at $x=0.5$ separates two
different states of the gas and the dimensionless initial states are
\begin{equation*}
(\rho,u,p)=\left\{\begin{aligned}
&(120, 0, 120/\gamma), 0<x<0.5,\\
&(1.2, 0, 1.2/\gamma),  0.5<x<1,
\end{aligned} \right.
\end{equation*}
where $Re=200$ and Prandtl number $Pr=0.73$.

\begin{figure}[!h]
\centering
\includegraphics[width=0.56\textwidth]{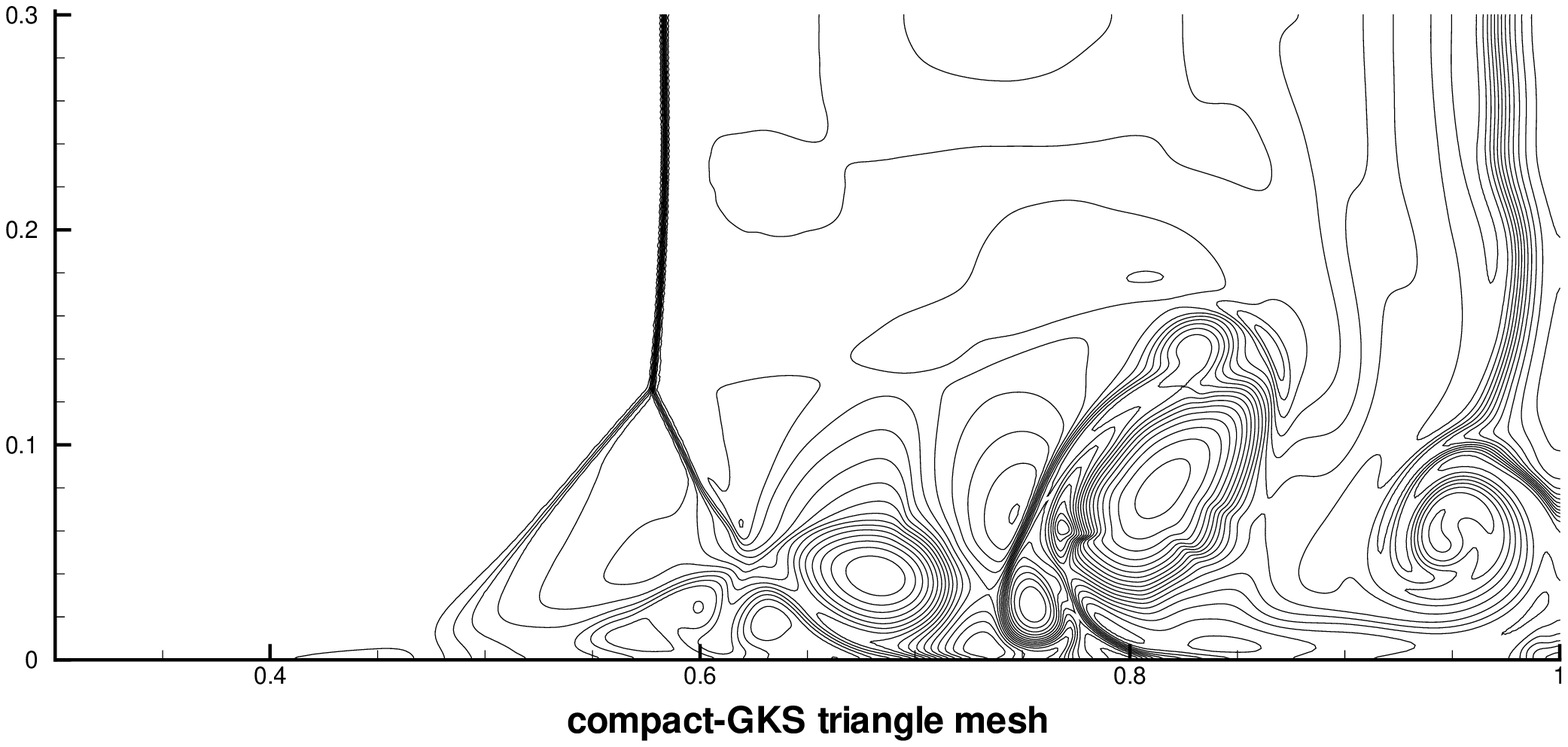}\\
\includegraphics[width=0.56\textwidth]{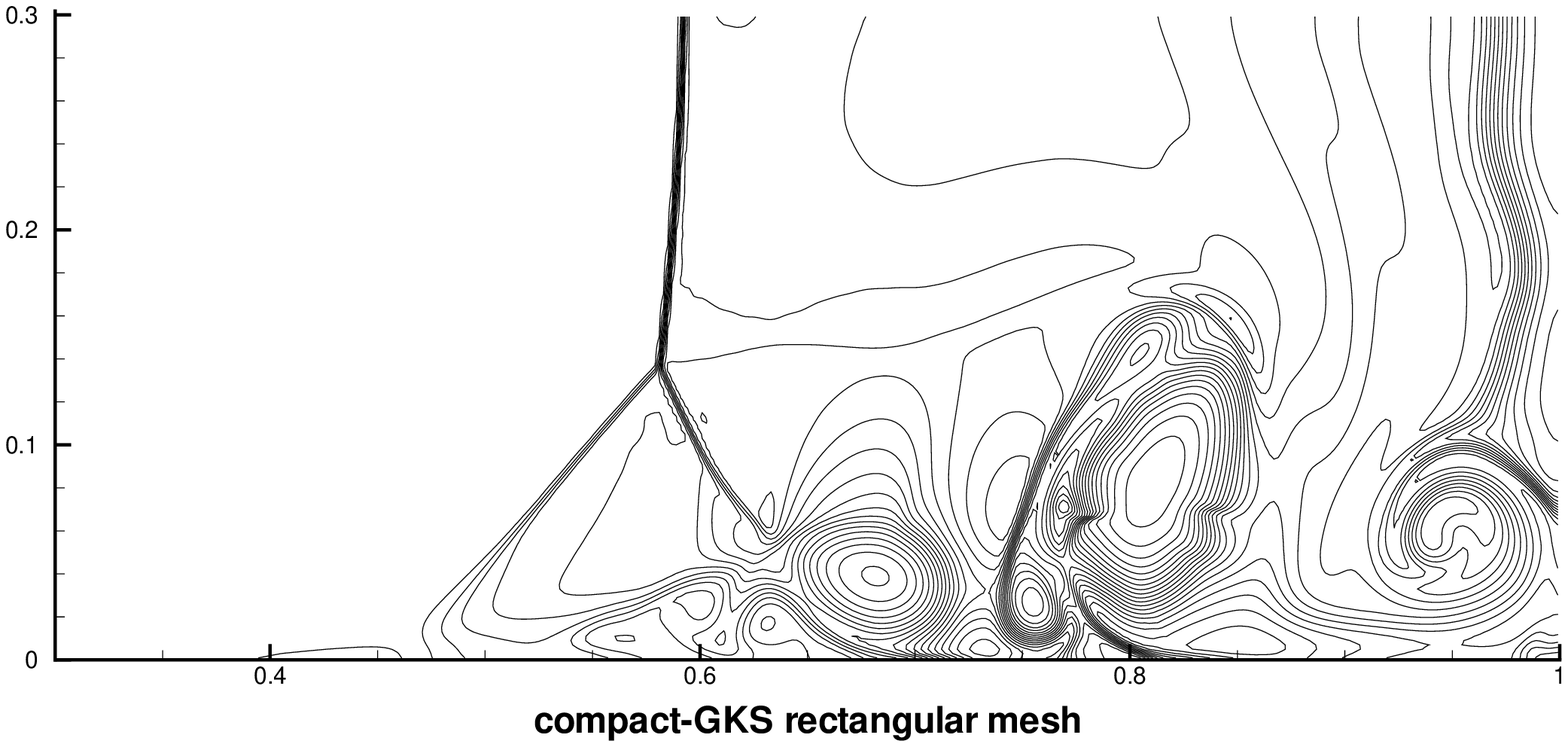}\\
\includegraphics[width=0.56\textwidth]{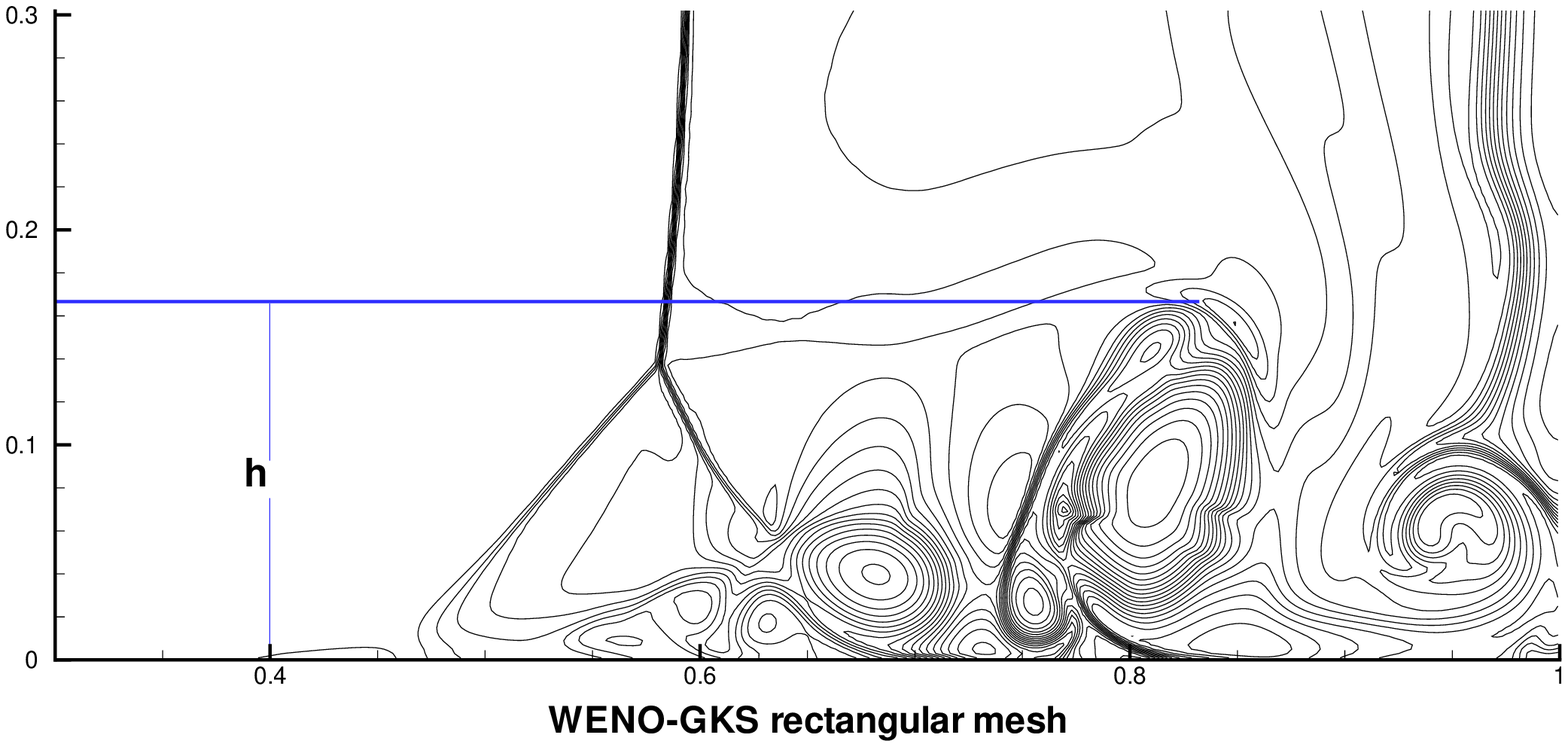}
\caption{\label{shock-boundary1} Reflected shock-boundary layer
interaction: the density distribution at $t=1$ with $Re=200$.}
\end{figure}

\begin{table}[!h]
\begin{center}
\def\temptablewidth{1\textwidth}
{\rule{\temptablewidth}{0.1pt}}
\begin{tabular*}{\temptablewidth}{@{\extracolsep{\fill}}cccccc}
Scheme  & AUSMPW+  &  M-AUSMPW+  & WENO-GKS & triangular & rectangular\\
\hline height & 0.163& 0.168& 0.165 & 0.164 & 0.166
\end{tabular*}
{\rule{\temptablewidth}{0.1pt}}\caption{\label{height} Comparison of
the height of primary vortex between gas kinetic schemes and
reference data \cite{Case-Kim} for the reflected shock-boundary
layer interaction.}
\end{center}
\end{table}

The membrane is removed at time zero and wave interaction occurs. A
shock wave, followed by a contact discontinuity, moves to the right
with a Mach number $Ma=2.37$, and reflects at the right end wall. After
the reflection, it interacts with the contact discontinuity. The
contact discontinuity and shock wave also interact with the horizontal
wall and create a thin boundary layer during their propagation. The
solution will develop complex two-dimensional
shock/shear/boundary-layer interactions. This case is tested in the
computational domain $[0, 1]\times[0, 0.5]$. A symmetrical condition
is used on the top boundary $x\in[0, 1], y=0.5$, and non-slip
boundary condition and adiabatic condition for temperature are
imposed at solid wall boundaries. The density distributions for the
compact scheme with the rectangular mesh with mesh size $\Delta
x=\Delta y=1/500$ and the triangular mesh with mesh size $h=1/500$ are
given in Fig.\ref{shock-boundary1}. As a reference, the density
distribution of the WENO-GKS with a structured mesh size $\Delta x=\Delta
y=1/500$ is also presented \cite{GKS-high3}. The current scheme can resolve the
complex flow structure. As shown in Table.\ref{height}, the height
of primary vortex predicted by the current scheme agrees well with
the reference data \cite{Case-Kim}. and the GKS-WENO results with the
same structure mesh.

\section{Conclusion}
In this paper, a third-order compact gas-kinetic scheme is proposed
on unstructured mesh for both inviscid and viscous flow simulations.
The merit of the gas-kinetic scheme is that due to a higher-order
gas evolution model the time-dependent solution of gas distribution
at a cell interface can provides both numerical fluxes and the
point-wise flow variables. Therefore, the scheme can be designed in
a compact way, where both the cell averaged and cell interface flow
variables can be used for the initial data reconstruction at the
beginning of next time level. With the inclusion of neighboring
cells only, a compact third-order gas-kinetic scheme is constructed,
where the weighted least-square method is used for the data
reconstruction on the unstructured mesh. In comparison with former
compact gas-kinetic scheme, the use of least-square procedure avoids
the difficulty in choosing different stencils. The systematic way of
including all weighted stencils makes the compact reconstruction
suitable for different kind of meshes. Different from other
higher-order schemes based on the Riemann solution, the current
method avoids the use of Gaussian points integration for the flux
transport along a cell interface and the multi-stage Runge-Kutta
time stepping technique. The compact scheme has been tested from
smooth viscous flow to the cases with strong discontinuities. The
numerical results confirm the accuracy and robustness of the current
third-order compact scheme.

\section*{Acknowledgement}
The work was supported by Hong Kong research grant council (620813,
16211014, 16207715).


\begin{thebibliography}{}
\bibitem{un-ENO} R. Abgrall  On essentially non-oscillatory schemes on
unstructured meshes: analysis and implementation. J. Comput. Phys.
144 (1994) 45-58.

\bibitem{k-exact-1} T.J. Barth, P.O. Frederichson,  Higher order solution of the
Euler equations on unstructured grids using quadratic
reconstruction. AIAA (1990) 90-0013.

\bibitem{BGK-1} P.L. Bhatnagar, E.P. Gross, M. Krook, A Model for
Collision Processes in Gases I: Small Amplitude Processes in Charged
and Neutral One-Component Systems, Phys. Rev. 94 (1954) 511-525.



\bibitem{BGK-3} S. Chapman, T.G. Cowling, The Mathematical theory of
Non-Uniform Gases, third edition, Cambridge University Press,
(1990).


\bibitem{Case-Daru} V. Daru, C. Tenaud, High order one-step monotonicity-preserving
schemes for unsteady compressible flow calculations, J. Comput.
Phys. 193 (2004) 563-594.



\bibitem{un-WENO3} M. Dumbser, M. K\"aser, V.A. Titarev, E.F. Toro. Quadrature-free
non-oscillatory finite volume schemes on unstructured meshes for
nonlinear hyperbolic systems, J. Comput. Phys. 226 (2007), 204-243.


\bibitem{un-WENO2} O. Friedrich, Weighted essentially non-oscillatory schemes for
the interpolation of mean values on unstructured grids, J. Comput.
Phys. 144 (1998) 194-212.




\bibitem{Case-Ghia} U. Ghia, K. N. Ghia, C.T Shin, High-Re solutions for incompressible
flow using the Navier-Stokes equations and a multigrid method, J.
Comput. Phys. 48 (1982) 387-411.

\bibitem{UGKS-Guo} Z.L. Guo, K. Xu, R.J. Wang,
Discrete unified gas kinetic scheme for all Knudsen number flows:
Low-speed isothermal case, Physical Review E, 88 (2013) 033305.



\bibitem{GKS-Jiang} J. Jiang, Y.H. Qian, Implicit gas-kinetic BGK scheme with
multigrid for 3D stationary transonic high-Reynolds number flows,
Computers $\&$ Fluids. 66 (2012) 21-28.

\bibitem{WENO2} G.S. Jiang, C. W. Shu, Efficient implementation of Weighted ENO
schemes, J. Comput. Phys. 126 (1996) 202-228.

\bibitem{Case-Kim} K.H. Kim, C. Kim, Accurate, efficient and monotonic numerical methods
for multi-dimensional compressible flows Part I: Spatial
discretization, J. Comput. Phys. 208 (2005) 527-569.

\bibitem{Shock-detection} L. Krivodonova, J. Xin, J.F. Remacle, N. Chevaugeond, J.E.
Flahertyd, Shock detection and limiting with discontinuous Galerkin
methods for hyperbolic conservation laws,  Applied Numerical
Mathematics 48 (2004) 323-338.

\bibitem{GKS-Kumar} G. Kumar, S.S. Girimaji, J. Kerimo,
WENO-enhanced gas-kinetic scheme for direct simulations of
compressible transition and turbulence, J. Comput. Phys. 234 (2013)
499-523.



\bibitem{GKS-high2} Q. Li, K. Xu, S. Fu,  A high-order gas-kinetic Navier-Stokes
flow solver, J. Comput. Phys. 229 (2010) 6715-6731.




\bibitem{GKS-high1} J. Luo, L.J. Xuan, and K. Xu,
Comparison of fifth-order WENO scheme and WENO-gas-kinetic scheme
for inviscid and viscous flow simulation, Commun. Comput. Phys.,  14
(2013) 599-620.


\bibitem{GKS-high3} J. Luo, K. Xu, A high-order multidimensional gas-kinetic scheme for hydrodynamic
equations, SCIENCE CHINA Technological Sciences, 56 (2013)
2370-2384.

\bibitem{UGKS-Luc} L. Mieussens, On the asymptotic preserving property of the unified gas-kinetic
scheme for the diffusion limit of linear kinetic models, J. Comput.
Phys. 253 (2013) 138-156.


\bibitem{k-exact-3} C. F. Ollivier-Gooch  Quasi-ENO schemes for
unstructured meshes based on unlimited data-dependent least-square
reconstruction,  J. Comput. Phys., 133 (1997) 6-17.

\bibitem{GKS-high4} L.Pan, K. Xu  A compact third-order gas-kinetic scheme for compressible Euler and
Navier-Stokes equations, Commun. Comput. Phys. 18 (2015) 985-1011.

\bibitem{Case-Pandolfi} M. Pandolfi, and D. D'Ambrosio,
Numerical Instabilities in Upwind Methods: Analysis and Cures for
the "Carbuncle" Phenomenon, J. Comput. Phys.  166 (2001) 271-301.



\bibitem{Case-Woodward} P. Woodward, P. Colella, The numerical simulation of two
dimensional fluids with strong shock, J. Comput. Phys. 54 (1984)
115-173.


\bibitem{GKS-Xu1} K. Xu, Direct Modeling for Computational Fluid Dynamics: Construction and Application
of Unified Gas-kinetic Schemes, World Scientific (2015).

\bibitem{GKS-Xu2} K. Xu, A gas-kinetic BGK scheme for the Navier-Stokes
equations and its connection with artificial dissipation and Godunov
method, J. Comput. Phys. 171 (2001) 289-335.


\bibitem{UGKS-Xu} K. Xu, J. Huang, A unified gas-kinetic scheme for continuum and
rarefied flows, J. Comput. Phys.  229 (2010) 7747-7764.
\end{thebibliography}
\end{document}